\documentclass[12pt,preprintnumbers]{article}
\usepackage{graphicx, amstext, amsmath, amssymb,color,amsthm}
\usepackage{enumerate}
\usepackage[font=scriptsize]{caption}
\usepackage{float}
\usepackage{hyperref}
\usepackage{fancyhdr}
\theoremstyle{plain}
\newtheorem{thm}{Theorem}

\newtheorem{prop}{Proposition}

\renewenvironment{proof}{{\bfseries Proof }} {\qed}

\title{Tail fitting for truncated and non-truncated
Pareto-type distributions}
\author{ Beirlant J.$^{a}$\footnote{Corresponding author: Jan Beirlant, KU Leuven, Dept of Mathematics and  LStat, Celestijnenlaan 200B, 3001 Heverlee, Belgium; Email: jan.beirlant@wis.kuleuven.be },  Fraga Alves, M.I.$^{b}$,  Gomes, M.I.$^{b}$,  \\
{$^a$ \fontsize{8pt}{11pt} \selectfont Department of Mathematics and Leuven Statistics Research Center, KU Leuven}
\\
{$^b$ \fontsize{8pt}{11pt} \selectfont Department of Statistics and Operations Research, University of Lisbon }
 }

\begin{document}

 \maketitle
\begin{abstract}
{\noindent Recently some papers, such as Aban, Meerschaert and Panorska (2006), Nuyts (2010) and Clark (2013),  have drawn attention to possible truncation in Pareto tail modelling.
Sometimes natural upper bounds exist that truncate the probability tail, such as the Maximum Possible Loss in insurance treaties. At other instances  ultimately at the largest  data, deviations from a Pareto tail behaviour become apparent.
This matter is especially important when extrapolation outside the sample is required.
Given that in practice one does not always know whether the distribution is truncated or not,  we consider estimators for extreme quantiles both under truncated and non-truncated Pareto-type distributions.
Hereby we make use of the estimator of the tail index for the truncated Pareto distribution first proposed in Aban {\it et al.} (2006).
We also propose a truncated Pareto QQ-plot and a formal test for truncation in order to help deciding between a truncated and a non-truncated case.
In this way we enlarge  the possibilities of extreme value modelling using Pareto tails, offering an alternative scenario by adding a truncation point $T$ that is large with respect to the available data.
 In the mathematical modelling we hence let $T \to \infty$ at different speeds compared to the limiting fraction ($k/n \to 0$) of data used in the extreme value estimation.
This work is motivated using practical examples from different fields of applications, simulation results, and some asymptotic results. }
\end{abstract}

\section{Introduction}

\begin{figure}[!ht]
    \begin{center}
   \includegraphics[width=\textwidth]{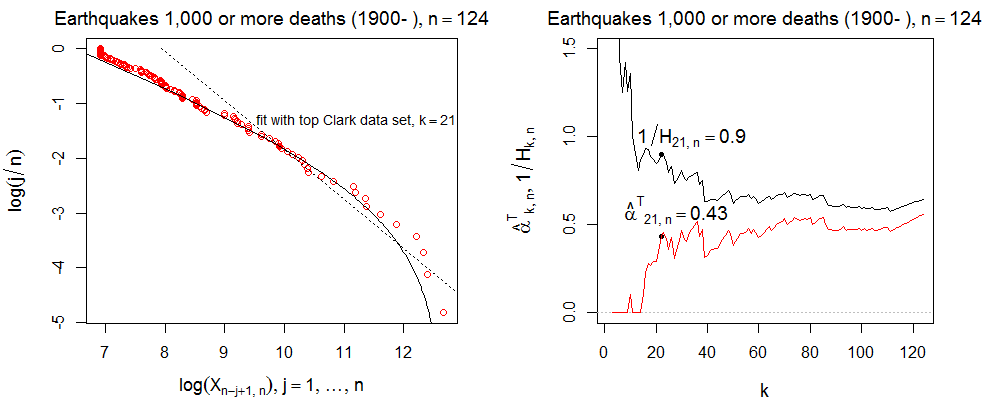}
  \caption {\small\scriptsize Earthquake fatalities data set. Left: Pa QQ-plot   with extrapolations anchored at $\log(X_{n-21,n}) $ based on a non-truncated Pa model in \eqref{powerlaw}  (dotted line) and a truncated Pareto model in \eqref{FtruncPar} (full line) as proposed in Clark (2013). Right: plot of $\hat{\alpha}^T_{k,n}$ and $1/H_{k,n}$ for $k=1,\ldots,n$.}
       \label{fig:1}
  \end{center}
   \end{figure}

 \noindent
Considering positive data, the Pareto (Pa) distribution is a simple and very popular model with power law probability tail.
 Using the notation from Aban {\it et al.} (2006), the right tail function (RTF)
 \begin{equation}
 \mathbb{P}(W>w) = \tau^{\alpha} w^{-\alpha} \mbox{ for } w \geq \tau >0 \mbox{ and } \alpha >0
 \label{powerlaw}
 \end{equation}
is considered as the standard example in the max domain of attraction of the Fr\'echet distribution. For instance, losses in property and casualty insurance often have a heavy right tail
behaviour making it appropriate for including large events in applications such as excess-of-loss pricing and enterprise risk management (ERM). There might be some practical problems with the use of
the Pa distribution and its generalization to the Pa-type model, because some probability mass can still be assigned to loss amounts that are unreasonable large or even physically impossible. In ERM this corresponds  to the concept of maximum possible loss. Here we will consider specific data sets on earthquake fatalities and forest fires. For other applications of natural truncation, such as probable maximum precipitation, see Aban   {\it et al.} (2006). These authors considered the upper-truncated Pareto distribution with RTF
\begin{equation}
\mathbb{P}(X>x) = \bar{F}_T (x) = \frac{\tau^{\alpha}(x^{-\alpha}-T^{-\alpha})}{1- (\tau/T)^{\alpha}}
\label{FtruncPar}
\end{equation}
for $0<\tau \leq x \leq T \leq \infty$, where $\tau < T$. Next to  the RTF, for any given distribution, we make use
 of the tail quantile function $U$ defined by $U(x) = Q(1-1/x)$ ($x>1$) where  $Q(1-p) := \inf \{x: F(x)\geq 1-p \}$ ($0<p<1$) denotes the upper quantile function corresponding to $F$. Henceforth, $\bar{F}_W$ and $U_W$, respectively $\bar{F}_T$ and $U_T$,  denote the RTF and the tail quantile function of the underlying Pa distribution, respectively of the Pa distribution truncated at $T$ from which the  $X$ data are observed with $X =_d W|W<T$. Note the following relations between the RTFs and the tail quantile functions:
\begin{eqnarray}
 \bar{F}_T (x) &=& \frac{\bar{F}_W(x)-\bar{F}_W(T)}{1-\bar{F}_W(T)}, \label{FWT}\\
 U_T (y) &=& U_W \left( y C_T \left[1 +yD_T \right]^{-1}\right)
 \label{UWTC}\\
 &= &U_W \left( {1 \over \bar{F}_W (T)}\left[1 +{1 \over yD_T} \right]^{-1}\right), \label{UWTD}
\end{eqnarray}
where $D_T = \bar{F}_W (T)/F_W(T)$ equals the odds ratio of the truncated probability mass under the untruncated Pa-type distribution $W$, and $C_T = 1/F_W(T)$.
\\\\
 Aban {\it et al.} (2006) derived the conditional maximum likelihood estimator (MLE) based on the $k+1$ ($0\leq k <n$) largest order statistics representing only the portion of the tail where the truncated Pareto (TPa) approximation holds. They showed that, with $X_{1,n} \leq \ldots \leq X_{n-k,n} \leq X_{n-k+1,n} \leq \ldots \leq X_{n,n}$ denoting the order statistics of an independent and identically  distributed sample of size $n$ from $X$, the MLE's obtained under this conditioning model are given by
$$\hat{T}= X_{n,n},\: \hat \tau  = k^{1/\hat{\alpha}_{k,n}^T}X_{n-k,n}  \left(n - (n-k)(X_{n-k,n}/X_{n,n})^{\hat{\alpha}_{k,n}^T} \right)^{-1/\hat{\alpha}_{k,n}^T} $$
where $\hat{\alpha}_{k,n}^T$ solves the equation
\begin{equation}
H_{k,n} = {1 \over \hat \alpha_{k,n}^T}
+ \frac{R_{k,n}^{\hat{\alpha}_{k,n}^T}\log R_{k,n}}{1-R_{k,n}^{\hat{\alpha}_{k,n}^T}},
\label{Hilltruncate}
\end{equation}
where $H_{k,n} =  {1\over k} \sum_{j=1}^k (\log X_{n-j+1,n}- \log X_{n-k,n})$ is the Hill (1975) statistic and $R_{k,n} = X_{n-k,n}/X_{n,n}$.
This estimator $1/\hat{\alpha}_{k,n}^T$ can be considered as an extension  of Hill's (1975) estimator to the case of a  TPa distribution, while $H_{k,n}$ was introduced as an estimator of $1/\alpha$ when $T=\infty$.
\\\\
Independently, Nuyts (2010) considered truncation of Pa distributions and obtained an adaptation of the Hill (1975) estimator through the estimation of $\mathbb{E}(\log W |W \in [L,R] )$ for some $0<L<R$ finite, taking $W$ to be the strict Pa in (\ref{powerlaw}). Then replacing $R$ by the truncation point $T$ and $L$ by some appropriate threshold $t=Q_T(1-(k+1)/(n+1))$, as commonly considered in extreme value methodology,  we obtain in the spirit of Nuyts (2010) that
\begin{equation}
\mathbb{E}(\log W - \log t |W \in [t,T] )= \frac{\int_{t}^T \log ({w \over t}) d \bar{F} (w) }{\int_{t}^T d\bar{F}(w)} = {1 \over \alpha} + \frac{({t \over T})^{\alpha}\log ({t \over T})}{1- ({t \over T})^{\alpha}}.
\label{nuyts1}
\end{equation}
Estimating $T$ by the maximum $X_{n,n}$ of the sample, $t$ by an order statistic $X_{n-k,n}$, and $\mathbb{E}(\log W - \log t |W \in [t,T] )$ by $H_{k,n}$,
 we find \eqref{Hilltruncate} again.
The estimator proposed in Nuyts (2010) differs slightly from $\hat{\alpha}_{k,n}^T$, but can be shown to be asymptotically equivalent to $\hat{\alpha}_{k,n}^T$. Nuyts (2010) also considered trimming the estimator by deleting the  $r-1$ ($r>1$) top data, leading to the  generalization of \eqref{Hilltruncate}
\begin{eqnarray*}
{1\over k-r+1} \sum_{j=r}^k (\log X_{n-j+1,n}- \log X_{n-k,n}) && \\
 && \hspace{-3.5cm} = {1 \over \hat \alpha^T_{r,k,n}}
+ \frac{X_{n-k,n}/X_{n-r+1,n})^{\hat{\alpha}^T_{r,k,n}}\log (X_{n-k,n}/X_{n-r+1,n})}{1-(X_{n-k,n}/X_{n-r+1,n})^{\hat{\alpha}^T_{r,k,n}}}.
\end{eqnarray*}
The estimator $\hat{\alpha}^T_{r,k,n}$ provides a way to make the Hill (1975) estimator robust against outliers, but it will be less efficient than the estimator $\hat{\alpha}_{k,n}^T$ without trimming. While robustness  under Pa models has received quite some attention in the literature (see for instance Hubert {\it et al.}, 2013, and the references therein), we confine ourselves to the case $r=1$.
\\\\
The solution of \eqref{Hilltruncate} can be approximated using  Newton-Raphson iteration on the equation
\[f\left(\frac 1\alpha\right):= H_{k,n} -{1 \over \alpha} - \frac{R_{k,n}^{\alpha} \log R_{k,n}}{1-R^{\alpha}_{k,n}} = 0 \]
to get
\begin{eqnarray}
{1 \over \hat{\alpha}_{k,n}^{(l+1)}} &= & {1 \over \hat{\alpha}_{k,n}^{(l)}}  \nonumber \\
 & & +
{ H_{k,n}-(\hat{\alpha}_{k,n}^{(l)})^{-1} - {R_{k,n}^{\hat{\alpha}_{k,n}^{(l)}}\log R_{k,n} \over 1- R_{k,n}^{\hat{\alpha}_{k,n}^{(l)}}} \over  1- (\hat\alpha_{k,n}^{(l)})^{2} {R_{k,n}^{\hat{\alpha}_{k,n}^{(l)}}\log^2 R_{k,n} \over (1- R_{k,n}^{\hat{\alpha}_{k,n}^{(l)}})^2 }} , \; l=0,1,...
\label{onestep}
\end{eqnarray}
where for instance Hill's estimator can serve as an initial approximation: $\: \hat{\alpha}^{(0)}=1/H_{k,n}$.
\\\\
The main purpose of this paper is to consider tail  estimation for truncated and non-truncated distributions in   case the Pa behaviour starts to set in from an intermediate threshold $t= Q_T(1-(k+1)/(n+1))$ on. In case of no truncation, this setting has been formalized mathematically by the concept of Pa-type distributions, defined by
\begin{equation}
\bar{F}_W (w) = w^{-\alpha} \ell_F (w), \; \alpha >0,
\label{paretotype}
\end{equation}
where $\ell_F$ is a slowly varying function at infinity, i.e. $\lim _{t \to \infty} \ell_F(ty)/\ell_F(t) = 1$ for every $y>0$. Hence, under this model, with $W/t$ denoting a peak over a threshold $t$ with $W>t$,
$$
\mathbb{P} ( W/t >y |W>t) \to y^{-\alpha} \mbox{ as } t \to \infty, \mbox{ for every } y>1.
$$
It is well-known that then also
\begin{equation}
U_W (y) = y^{1/\alpha} \ell_U (y), \; y >1,
\label{paretotype2}
\end{equation}
where $\ell_U$ is again a slowly varying function.
\\
In extreme value statistics the parameter $\xi:=1/\alpha$  is referred to as the extreme value index (EVI). The EVI $\xi$ is the shape parameter in the generalized extreme value distribution
\begin{equation*}
G_{\xi}(x) = \exp \left( - (1+\xi x)^{-1/\xi}\right), \mbox{ for } 1+\xi x>0.
\end{equation*}
This class of  distributions  is the set of the unique non-degenerate limit distributions of a sequence of maximum values, linearly normalized. In case $\xi >0$ the class of distributions for which the maxima are attracted to $G_{\xi}$ corresponds to the Pa-type distributions in (\ref{paretotype}). Note that for a given $T$ fixed truncated Pa models are known to exhibit an EVI $\xi=-1$, see for instance Figure 2.8 in Beirlant {\it et al.} (2004).
\\\\
To illustrate the practical importance of the present setting, we consider the data set  containing fatalities due to large earthquakes as published by the U.S. Geological Survey on http://earthquake. usgs.gov/earthquakes/world/,  which were also used in Clark (2013).
It contains the estimated number of deaths for the 124 events between 1900 and 2011 with at least 1000 deaths.
In Figure \ref{fig:1} the Pa QQ-plot (or log-log plot)
\begin{equation}
\left(\log X_{n-j+1,n}, \log (j/n)  \right), \;\; j= 1,\ldots, n,
\label{PaQQ}
\end{equation}
 is given. It exhibits a linear or Pa pattern for a large section of the data, while the final end of the plot is curved.  The plot of $1/H_{k,n}$ as the classical estimates of $\alpha$ based on an untruncated Pa distribution indeed bends upwards towards smaller values of $k$. This indicates that the unbounded Pareto pattern could be violated in this example. On this plot the extrapolations using a Pa distribution \eqref{powerlaw} respectively the TPa model (\ref{FtruncPar})  (with the linear respectively the curved extrapolation) are plotted based on the largest 21 data points as it was proposed in Clark (2013) such that $\hat{\alpha}^T_{21,n} = 0.43$ and $1/H_{21,n}=0.90$.
In section 3 tail extrapolation based on  the TPa model will indeed appear to be appropriate in this case.

   \vspace{0.2cm} \noindent
As suggested in Clark (2013) $\alpha$ could also be taken to be zero or negative. For instance formally setting $\alpha =-1$ in \eqref{FtruncPar}, one obtains the tail of a uniform type distribution. Finite tail distributions following \eqref{FtruncPar} with negative value of $\alpha$ show a fast rate of convergence to $T$ for small $p$  and when $T$ is a big number, due to the presence of $D_T$ in this expression.  In the applications we have in mind here, as we allow $T$ to be large,  convergence to $T$ is slow and hence a positive value of $\alpha$ is  appropriate. Closely related to this paper, Chakrabarty and  Samorodnitsky (2012) considered a different truncation model, whereby a heavy tailed random variable is truncated at a high value $M \to \infty$ added by exponentially tailed random variables $R$, and where the probability mass under $W$ behind $M$ is set at $M+R$. These authors show the consistency of Hill's estimator and provide a test for this kind of truncation. They did not consider the estimation of extreme quantiles however.
\\\\
In the next sections we specify the truncation models and provide estimators for extreme quantiles and for $T$. In section 4 we consider the problem of deciding between a Pa-type  case in  (\ref{paretotype}) and a TPa-type case in (\ref{FT}), or, between light and rough truncation.  To this end we construct a TPa QQ-plot, provide a new significance test and compare it with the test proposed in Aban {\it et al.} (2006).
In section 5 we study the finite sample behaviour of the proposed estimators using simulations and practical examples. The asymptotic properties of $\hat{\alpha}^T_{k,n}$ and the extreme quantile estimators under rough and light truncation are discussed in a final section. Proofs will be deferred to the supplementary material.

\section{Rough and light truncation of Pareto-type distributions}

Truncation of a Pa-type distribution at a value $T$ necessarily requires $t<T \to \infty$ and
\begin{eqnarray*}
\mathbb{P}(X/t >y|X>t) &=& \mathbb{P}(W/t >y | t<W<T) \\ &= & \frac{(yt)^{-\alpha}\ell_F (yt)-T^{-\alpha}\ell_F (T) }{t^{-\alpha}\ell_F (t)-T^{-\alpha}\ell_F (T)}   \\
&=&  \frac{y^{-\alpha}{\ell_F (yt) \over \ell_F(t)} -\left({T \over t}\right)^{-\alpha}{\ell_F (T) \over \ell_F(t)} }{1-\left({T \over t}\right)^{-\alpha}{\ell_F (T) \over \ell_F(t)}}.
\end{eqnarray*}
One can now consider two cases as $t,T \to \infty$:
\begin{itemize}
\item {\it Rough truncation}  when $T/t \to \beta > 1$, and by Karamata's uniform convergence theorem for regularly varying functions (Seneta, 1976),
\begin{equation}
\mathbb{P}(X/t >y|X>t) \to {y^{-\alpha} - \beta^{-\alpha} \over 1- \beta^{-\alpha}}, \;\; 1< y < \beta.
\label{FT}
\end{equation}
This corresponds to situations where  the deviation from the Pa behaviour due to truncation at a high value will be visible in the data, and an adaptation of the
classical Pa tail extrapolation methods appears appropriate.
\item {\it Light truncation} when $T/t \to \infty$
\begin{equation}
\mathbb{P}(X/t > y |X>t)  \to y^{-\alpha}, \;\; y>1,
\label{LT}
\end{equation}
and hardly any truncation is visible in the data, and the Pa-type model \eqref{paretotype} without truncation and the existing extreme value methods for Pa-type tails are appropriate.
\end{itemize}

\vspace{0.2cm}\noindent
When $T \to \infty$ we have from\eqref{UWTD} that
\begin{eqnarray*}
{U_T (y) \over U_W (1/\bar{F}_W (T) )} &=&
{U_W \left( {1 \over \bar{F}_W (T)}\left[1 +{1 \over yD_T} \right]^{-1}\right) \over  U_W (1/\bar{F}_W (T) )} \\
 &=& \left[1 +{1 \over yD_T}\right]^{-1/\alpha}\zeta_{y,T}
 \end{eqnarray*}
 where
\begin{equation}
\zeta_{y,T}:=\frac{ \ell_U \left({1 \over \bar{F}_W (T)}\left[1 +{1 \over yD_T} \right]^{-1} \right)}{\ell_U \left({1 \over \bar{F}_W (T)}\right)}.
\label{zeta}
\end{equation}
Since $U_W (1/\bar{F}_W (T) )=T$ and $1/\bar{F}_W (T) \to \infty$, we now obtain that if $yD_T$ is bounded away from 0
\begin{equation}
Q_T (1-1/y)= U_T(y)=  T \left( 1+ {1 \over yD_T}\right)^{-1/\alpha} \zeta_{y,T},
\label{QT}
\end{equation}
where $\zeta_{y,T} \to 1$ as $T \to \infty$ using Karamata's uniform convergence theorem. Note that  $\zeta_{y,T} = 1$ in case $W$ satisfies \eqref{powerlaw}. \\
In case $yD_T \to 0$ it follows similarly from \eqref{UWTC} that
\begin{equation}
Q_T (1-1/y)= U_T(y)= (yC_T)^{1/\alpha}[1+ yD_T]^{-1/\alpha}\ell_U (y C_T [1+ yD_T]^{-1}).
\label{QTb}
\end{equation}
Hence, with $(n+1)/(k+1)$ playing the role of $y$, and choosing $t =Q_T (1-(k+1)/(n+1)) = U_T((n+1)/(k+1))$, the conditions of rough and light truncation are rephrased as follows, in accordance with \eqref{QT} and \eqref{QTb}:
\begin{itemize}
\item {\it Rough truncation} when $$T/ U_T((n+1)/(k+1))  \to \beta > 1, \mbox{ or } {k+1 \over (n+1)D_T} \to \kappa :=\beta^{\alpha}-1,$$
and
\begin{equation}
{U_T (u{n+1 \over k+1}) \over U_T({n+1 \over k+1})} \to  \left( 1+\kappa \over 1 +\kappa /u \right)^{1/\alpha},\;\; u>1;
\label{UT}
\end{equation}
\item {\it Light truncation} when
$$T/ U_T((n+1)/(k+1))  \to \infty \mbox{ or } {(n+1)D_T \over k+1} \to  0$$
and
$$
{U_T (u{n+1 \over k+1}) \over U_T({n+1 \over k+1})} \to  u ^{1/\alpha}, \;\; u>1.
$$
\end{itemize}

\vspace{0.2cm}\noindent
The estimation of $\alpha$ and $k/(nD_T)$ (or $\kappa$, $\beta$) will now constitute important steps in order to arrive at estimators of extreme quantiles $Q_T (1-p)$ and $T$.
Given the fact that our model is only defined choosing $k,n,T \to \infty$, $k/n \to 0$ jointly with the corresponding conditions for rough and light truncation, the underlying model depends on $n$ and a triangular array formulation  $X_{n1}, \ldots, X_{nn}$ of the observations should be used in order to emphasize the  nature of the model. However, in statistical procedures as presented here, when a single sample is given, the notation $X_1,\ldots,X_n$ is more natural and will be used throughout.
\\\\

\section{Estimation of extreme quantiles}

From (\ref{QT}) it is clear that the estimation of $D_T$ is an intermediate step in important estimation problems following the estimation of $\alpha$, namely of extreme quantiles and of the endpoint $T$.
From (\ref{QT})
\begin{equation}
\left( {Q_T(1-{k+1 \over n+1}) \over Q_T(1-{1 \over n+1})} \right)^{\alpha} \approx { 1+ {1 \over (n+1) D_T} \over  1+ {k+1 \over (n+1)D_T}} = {1 \over k+1} \left( {1+ (n+1)D_T  \over  1+ {(n+1)D_T \over k+1}}\right).
\label{ratio1}
\end{equation}
Motivated by (\ref{ratio1}) and estimating $ Q_T(1-(k+1)/(n+1)) / Q_T(1-1/(n+1))$ by $R_{k,n}$, we propose
\begin{equation}
\hat{D}_T   := \hat{D}_{T,k,n} =  {k+1 \over n+1}{ R_{k,n}^{\hat{\alpha}^T_{k,n}} - {1 \over k+1}
\over 1-R_{k,n}^{\hat{\alpha}^T_{k,n}}}
\label{DT}
\end{equation}
 as an estimator for $D_T$ in case of truncated and non-truncated Pa-type distributions.
In practice we will make use of the admissible estimator
$$
\hat{D}^{(0)}_T :=  \max \left\{ \hat{D}_T,0\right\}.
$$
\\
In case $D_T >0$, in order to construct estimators of $T$ and extreme quantiles $q_p=Q_T(1-p)$,
as in (\ref{ratio1}) we find that
\begin{equation}
\left( {Q_T(1-p) \over Q_T(1-{k+1 \over n+1})} \right)^{\alpha}
\approx { 1+ {k+1 \over (n+1)D_T} \over  1+ {p \over D_T} } =  { D_T+ {k+1 \over n+1} \over  D_T+ p  } .
\label{ratio2}
\end{equation}
Then  taking logarithms on both side of (\ref{ratio2}) and estimating  $Q_T(1-(k+1)/(n+1))$ by $X_{n-k,n}$ we find an estimator   $\hat{q}^T_p := \hat{q}^T_{p,k,n}$ of $q_p$:
\begin{equation}
\log \hat{q}^T_{p,k,n} = \log X_{n-k,n} + {1 \over \hat{\alpha}^T_{k,n}}
\log \left( {\hat{D}_T+ {k+1 \over n+1} \over \hat{D}_T + p} \right).
\label{Qest1}
\end{equation}
Note that $\hat{q}^T_p$ can also be rewritten as
\begin{eqnarray}
\log \hat{q}^T_{p,k,n} &=& \log X_{n-k,n} + {1 \over \hat{\alpha}^T_{k,n}}
\log \left( {1+ {k+1 \over (n+1)\hat{D}_{T}} \over 1+ {p \over \hat{D}_{T}}} \right),
\label{Qest2} \\
 \hat{q}^T_{p,k,n} &=& X_{n-k,n}\left( {k+1 \over (n+1)p} \right)^{1/\hat{\alpha}^T_{k,n}}
\left( {1+ { (n+1)\hat{D}_{T} \over k+1} \over 1+ { \hat{D}_{T}\over p}} \right)^{1/\hat{\alpha}^T_{k,n}}.
\label{QTinf}
\end{eqnarray}
An estimator of $T$ follows from letting $p \to 0$,
$$\log \hat{T}_{k,n}=\max \left\{\log \hat{q}^T_{0,k,n} , \log X_{n,n}\right\},$$
taking the maximum with  $\log X_{n,n}$ in order for this endpoint estimator to be admissible.\\
Equation (\ref{QTinf}) for $\hat{q}^T_{p,k,n}$ constitutes an adaptation to the TPa case of the Weissman (1978) estimator
\begin{equation}
\hat{q}_{p,k,n}^{W} = X_{n-k,n}\left( {k+1 \over (n+1)p} \right)^{H_{k,n}}
\label{W}
\end{equation}
which is valid  under \eqref{paretotype}. Expression \eqref{QTinf} is more adapted to the case of light truncation or $(n+1)D_T/(k+1) \to 0$. Version (\ref{Qest2}) can be linked to the case of rough truncation or $(k+1)/((n+1)D_T) \to \kappa$.
Version (\ref{Qest1})  can be applied in all cases. However in section 5, in case of light truncation $\hat{q}^T_{p,k,n}$ will be shown to be consistent only when  $(np_n)^{-1} \to 0$. When this condition is not satisfied under light truncation we propose to use
\begin{equation}
\log \hat{q}_{p,k,n}^{\infty} = \log X_{n-k,n}+ {1 \over\hat{\alpha}^T_{k,n}}\log \left( {k+1 \over (n+1)p} \right).
\label{QTWT}
\end{equation}
Note that such alternative expressions do not exist for the estimation of the endpoint $T$ as in case $D_T=0$ no finite endpoint exists.

\section{Goodness-of-fit and testing for truncated Pareto-type distributions}

From the preceeding sections the need for goodness-of-fit and test for truncated Pa distributions became apparent.
Based on a chosen value  $\hat{D}_{T,k^*,n}$ for particular $k^*$, we propose the TPa QQ-plot to verify the validity of (\ref{QT}):
\begin{equation}
\left(\log X_{n-j+1,n}, \log \left(\hat{D}_{T,k^*,n}+ j/n \right)  \right), \;\; j= 1,\ldots, n.
\label{TPQQ}
\end{equation}
Note that when $T=\infty$ or $D_T=0$ the TPa QQ-plot  agrees with the classical Pareto QQ-plot \eqref{PaQQ}.
Under (\ref{QT}) an ultimately linear pattern should be observed to the right of some anchor point, i.e. at the points with indices $j=1,\ldots,k$ for some $1< k <n$. From this, we propose to choose the value of $k^*$ in practice
as the value that maximizes the correlation between $\log X_{n-j+1,n}$ and $ \log \left(\hat{D}_{T,k^*,n}+ j/n \right) $ for $j=1,\ldots,k^*$ and $k^*>10$. This choice can be improved in future work since the covariance structure of the deviations of the points on the TPa QQ-plot from the reference line are neither independent nor identically distributed. This issue was addressed for the Pa QQ-plot in Beirlant {\it et al.} (1996) and Aban and Meerschaert (2004) and should be considered in the truncated  case too.
\\\\
Aban {\it et al.} (2006) already proposed a test for $H^{(1)}_0: T = \infty$ versus $H^{(1)}_1: T<\infty$ under the strict Pa and TPa models, rejecting $H_0^{(1)}$ at asymptotic level $q \in (0,1)$ when
\begin{equation}
X_{n,n} < \left( {nA \over -\log q}\right)^{1/\alpha}
\label{Atest}
\end{equation}
for some $1<k<n$ and where $A=\tau^{\alpha}$ in \eqref{powerlaw}. In \eqref{Atest},  $\alpha$ is estimated by the maximum likelihood estimator $1/H_{k,n}$ under $H_0^{(1)}$ based on the Hill (1975) estimator
while
\begin{equation}
\hat{A}_{k} =  {k \over n} \left(  X_{n-k,n}\right)^{1/H_{k,n}}.
\label{1overHill}
\end{equation}
Note that the rejection rule \eqref{Atest} can be rewritten as
\begin{equation}
T_{A,k,n} := kR_{k,n}^{1/H_{k,n}} >\log {1 \over q},
\label{Atest2}
\end{equation}
and the p-value is given by $\exp (-k R_{k,n}^{1/H_{k,n}})$.
\\\\
Here we  also consider the problem of testing light versus rough truncation, i.e.
$$H^{(2)}_{0}:X \mbox{ satisfies }  \eqref{LT} \mbox{ {\it versus} } H^{(2)}_{1}:X \mbox{ satisfies }  \eqref{FT}.$$
In the next sections we inspect the finite sample and asymptotic  properties of the test \eqref{Atest} under $H^{(2)}_0$ and $H^{(2)}_1$.
\\
We propose also a different test  rejecting $H^{(2)}_0$ when an appropriate estimator of $(n+1)D_T/(k+1)$ is significantly different from 0.  Here we construct such an estimator generalizing $R_{k,n}^{\alpha}$ with an average of ratios $(X_{n-k,n}/X_{n-j+1,n})^{\alpha}, \; j=1,\ldots,k$, which then possesses an asymptotic normal distribution under the null hypothesis. Observe that  under \eqref{QT} as $k \to \infty$
\begin{eqnarray*}
MR_{k,n}={1 \over k} \sum_{j=1}^k \left( {Q_T (1- {k+1 \over n+1})  \over Q_T (1-{j \over n+1}) }\right)^{\alpha} &\approx &
{1 \over k} \sum_{j=1}^k {1+{j \over k+1}{k+1 \over (n+1)D_T} \over 1+{k+1 \over (n+1)D_T}} \\
&\approx & {1+{1 \over 2 }{k+1 \over (n+1)D_T}  \over 1+{k+1 \over (n+1)D_T} }.
\end{eqnarray*}
Estimating  $MR_{k,n}$ by
$$
E_{k,n}(\alpha) = {1 \over k} \sum_{j=1}^k \left( {X_{n-k,n} \over X_{n-j+1,n} }\right)^{\alpha},
$$
leads now to
\begin{equation}
L_{k,n}(\hat\alpha) := {E_{k,n}(\hat\alpha) -{1 \over 2} \over 1- E_{k,n}(\hat\alpha)}
\label{L}
\end{equation}
as an estimator of $(n+1)D_T/(k+1)$, with $\hat{\alpha}$ an appropriate estimator of $\alpha$.
Under $H^{(2)}_0$ the reciprocal of the Hill (1975) estimator $1/H_{k,n}$ is an appropriate estimator of $\alpha$.   Moreover, in the final section it will be stated that under some regularity assumptions on the underlying Pa-type distribution, we have under $H^{(2)}_0$ for $k,n \to \infty$ and $k/n\to 0$, that $ \sqrt{k} L_{k,n}(1/H_{k,n})$ is asymptotically normal with mean 0 and variance 1/12.
Moreover,  it is then also shown that under rough truncation
as $k,n,T \to \infty$, $k/n \to 0$ and $nD_T/k \to  1/\kappa >0$
$$
L_{k,n}(H^{-1}_{k,n}) \to_p (1+\kappa)^{-{\kappa \over \kappa - \log (1+\kappa)}}{\kappa - \log (1+\kappa) \over \kappa (2\kappa - \log (1+\kappa) )}
\left( (1+\kappa)^{{2\kappa - \log (1+\kappa) \over \kappa - \log (1+\kappa)}} -1 \right) <0,
$$
so that an asymptotic test based on $L_{k,n}(1/H_{k,n})$ rejects $H^{(2)}_0$ on asymptotic level $q$ when
\begin{equation}
T_{B,k,n}:= \sqrt{12 \, k} L_{k,n}(1/H_{k,n}) < -z_q
\label{Btest}
\end{equation}
with $P(\mathcal{N}(0,1) > z_q)=q$. The p-value is then given by $\Phi (\sqrt{12 \, k} L_{k,n}(1/H_{k,n}))$.
Both tests \eqref{Atest} and \eqref{Btest} will further be compared below.

\section{Practical examples and simulations}

\begin{figure}[!ht]
    \begin{center}
      \includegraphics[width=\textwidth]{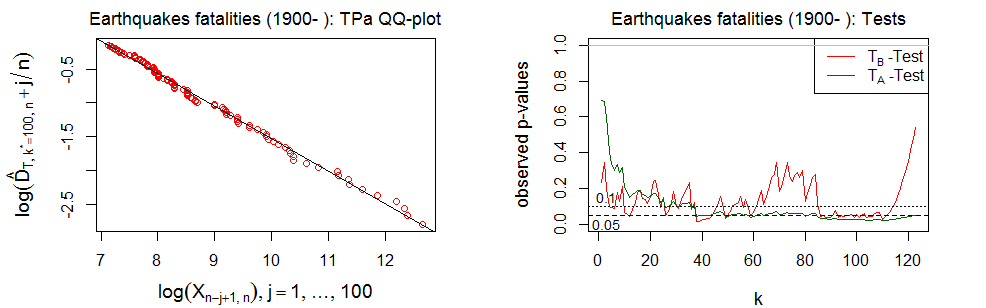}
\includegraphics[width=\textwidth]{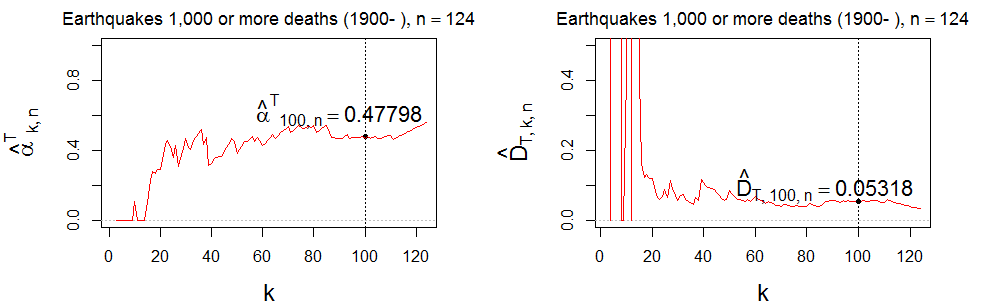}
  \includegraphics[width=\textwidth]{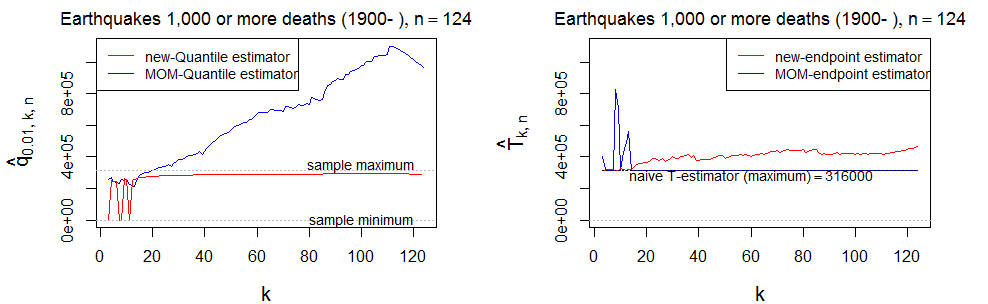}
\caption {\small\scriptsize Earthquake fatalities data set. Top:  TPa QQ-plot (left) for the earthquake fatalities data set using  $k^*=100$; plot of p-values based on $T_A$ and $T_B$. Middle: plots of Pareto index $\hat{\alpha}^T_{k,n}$ (left)
and  odds ratio $\hat{D}_{T,k,n}$ estimates ($k=1,\ldots,124$) (right) marking the values at $k=100$. Bottom: quantile estimates $\hat{q}^T_{0.01,k,n}$ (left)  and endpoint estimates $\hat{T}_{k,n}$  (right)  contrasted with the method of moments quantile and endpoint estimators,  in \eqref{quantile-MOM} and \eqref{endpoint-MOM}, respectively.}
       \label{fig:1b}
  \end{center}
   \end{figure}

First we retake the data set containing fatalities due to large earthquakes from Figure 1.
In Figure \ref{fig:1b} (\emph{middle}) the estimates $\hat{\alpha}^T_{k,n}$ and $\hat{D}_{T,k,n}$ are plotted against $k=1,\ldots,n$. Here we have chosen $k^* = 100$ as a typical value where both plots are horizontal in $k$.
The TPa QQ-plot in (\ref{TPQQ}) is given in Figure \ref{fig:1b} (\emph{top left}), using the above mentioned value  $k^*=100$. Next, the p-values as a function of $k$, both for the $T_A$ and $T_B$ tests are given (\emph{top right}).
Finally in Figure \ref{fig:1b} (\emph{bottom})  the estimates of the extreme quantile $q_{0.01}$ using (\ref{Qest1}) and the endpoint $T$ (obtained by letting $p \to 0$) are presented as a function of $k$. They are contrasted with the values obtained by the classical method of moment estimates as introduced in  Dekkers {\it et al.} (1989) illustrating the slow convergence of the classical extreme value methods in the TPa-type model we study here.
For any real EVI, the classical moment $\xi$-estimator is  defined by
\begin{equation}\label{MOM}
\hat \xi_{n,k}^{MOM} := M_{n,k}^{(1)}+\hat \xi_{n,k}^{-}, \quad \xi_{n,k}^{-} :=1-\frac{1}{2} \left[1-\left(M_{n,k}^{(1)}\right)^2/M_{n,k}^{(2)}\right]^{-1},
 \end{equation}
with $
M_{n,k}^{(j)}:=\frac{1}{k}\sum_{i=0}^{k-1}\ln^j \left( X_{n-i,n}/X_{n-k,n}\right)$, $j=1,2$, which constitutes a consistent estimator for $\xi \in \mathbb{R}$. The Hill estimator is $M_{n,k}^{(1)}=H_{k,n}$.
\\
The MOM-estimators for high quantiles and right endpoint, based on the moment estimator $\hat \xi_{n,k}^{MOM}$,  are defined by
(see de Haan and Ferreira, 2006, $\S 4.3.2$, for details)
 \begin{equation}\label{quantile-MOM}
\hat q_p^{MOM}:=X_{n-k,n}+X_{n-k,n}M_{n,k}^{(1)}\left(1-\hat \xi_{n,k}^{-}\right)\frac{\left(\frac{k}{np}\right)^{\hat \xi_{n,k}^{MOM}} -1}{\hat \xi_{n,k}^{MOM}}
 \end{equation}
 and
 \begin{equation}\label{endpoint-MOM}
\hat T^{MOM}:=\max\left(\, \,\hat T^{(M)} ,  X_{n,n}  \,\right),\,\, \hat T^{(M)}:= X_{n-k,n}-\tfrac{X_{n-k,n}M_{n,k}^{(1)}(1-\hat \xi_{n,k}^{-})}{\hat \xi_{n,k}^{MOM}}.
 \end{equation}
Notice that in \eqref{endpoint-MOM}  $\hat T^{MOM}$ corresponds to the admissible version of the moment endpoint estimator $\hat T^{(M)}$, since the  latter can return values below the sample maximum.
\\\\
Concerning the high quantile estimation, the chosen value  $p=0.01$ is directly related to the modest sample size here of $n=124$. The quantile estimates $\hat{q}^T_{0.01,k,n}$  reveal a stable pattern on $k$, in Figure \ref{fig:1b} (\emph{bottom left}).
\\While on the basis of Figure \ref{fig:1b}, $T_A$ is only border significant for small values of $k$, the TPa-type model with a truncation point $T$ around 400,000 deaths offers a convincing fit, and leads to a useful estimator for extreme quantiles.

\begin{figure}[!ht]
    \begin{center}
  \includegraphics[width=\textwidth]{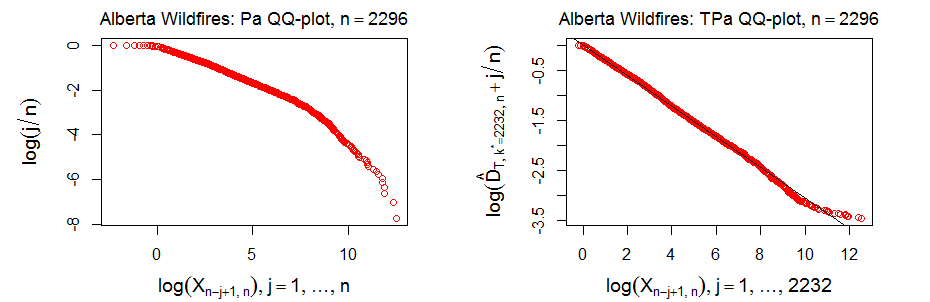}
  \includegraphics[width=\textwidth]{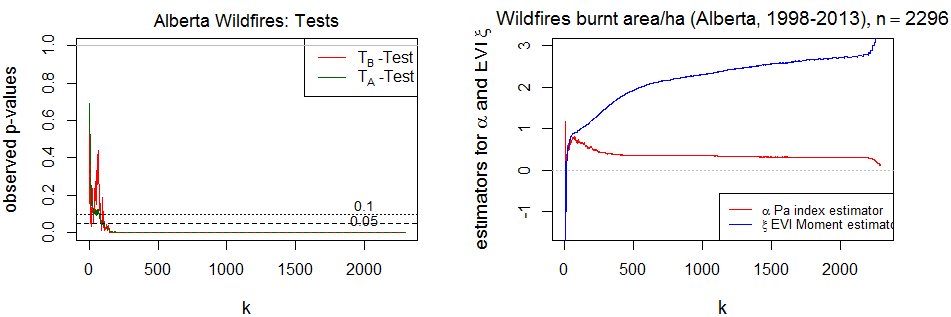}
  \includegraphics[width=\linewidth]{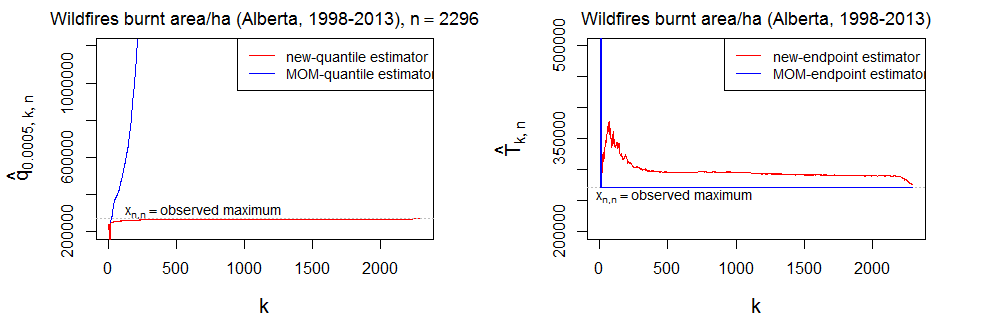}
 \caption{\small\scriptsize Alberta wildfires. Top:  Pa QQ-plot  (left); TPa QQ-plot (right) using  the top $k^*=2232$ data. Middle: plots of p-values for $T_A$ and $T_B$ tests (left), and  Pa index $\hat{\alpha}^T_{k,n}$ (right) estimates contrasted with the method of moments estimates in \eqref{MOM}. Bottom: quantile estimates $\hat{q}^T_{0.0005,k,n}$  (left)  and endpoint estimates $\hat{T}_{k,n}$ (right)  contrasted with the method of moments quantile and endpoint estimators,  in \eqref{quantile-MOM} and \eqref{endpoint-MOM}, respectively.}
       \label{fig:2a}
  \end{center}
   \end{figure}

\begin{figure}[!ht]
    \begin{center}
  \includegraphics[width=0.48\textwidth]{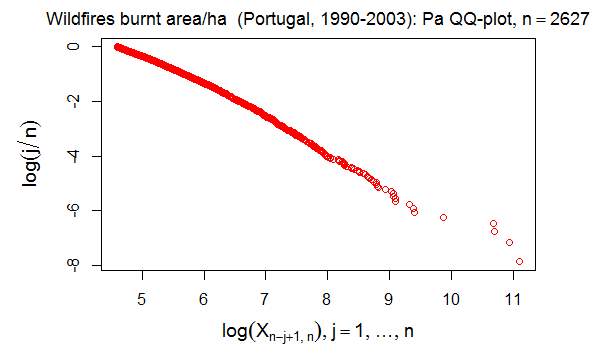}
   \includegraphics[width=\linewidth]{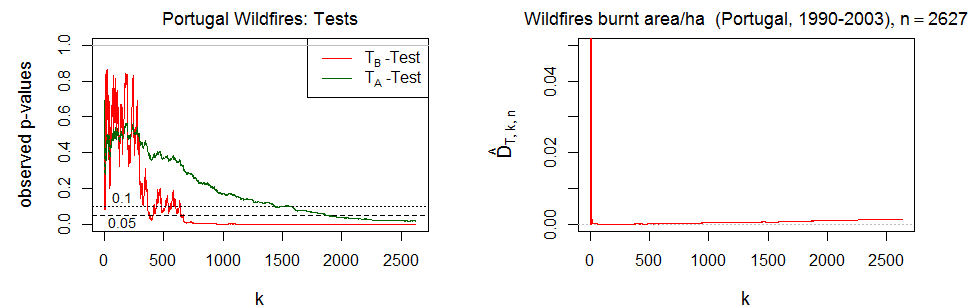}
  \includegraphics[width=\linewidth]{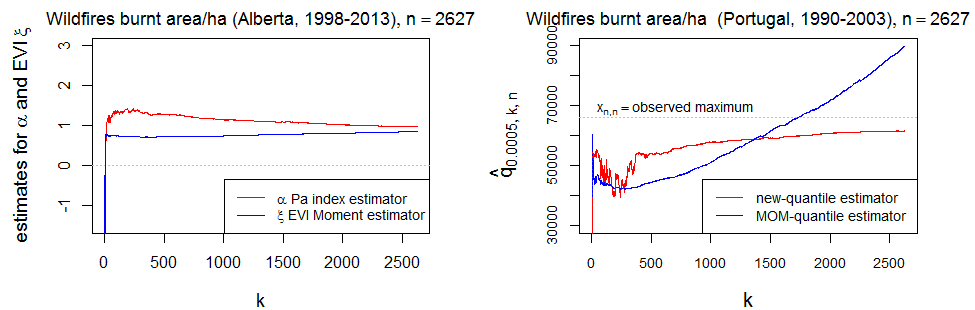}
\caption{\small\scriptsize Portugal wildfires. Top:  Pa QQ-plot. Middle: plots of p-values for $T_A$ and $T_B$ tests (left) , and  odds ratio $\hat{D}_{T,k,n}$ estimates  (right). Bottom: Pa index $\hat{\alpha}^T_{k,n}$ estimates (left) and quantile estimates $\hat{q}^T_{0.0005,k,n}$  (right)  contrasted with the method of moments index and quantile estimators,  in \eqref{MOM} and \eqref{quantile-MOM}, respectively.}
       \label{fig:2b}
  \end{center}
   \end{figure}

\vspace{0.3cm} \noindent
Another application can be found in statistical modelling of size distributions of forest fires. Power law distributions have appeared in literature for modelling such sizes, while Reed and McKelvey (2002) provided evidence that in some circumstances this is too simple to describe such distributions over their full range. We consider here two data sets with sizes of wild fires from Alberta (Canada) from 1998-2013 ($n=2296$) that can be found on
http://wildfire.alberta.ca/wildfire-maps/historical-wildfire-information/spatial-wildfire-data.aspx, together with  a data set, recently considered in Gomes \emph{et al. }(2012) and in Brilhante \emph{et al.} (2013), containing the number of hectares, exceeding 100 ha, burnt during wildfires recorded in Portugal from 1990 till 2003 ($n=2627$).
Both tails are analyzed in Figure \ref{fig:2a} for the Alberta data set and in  Figure \ref{fig:2b} for the Portuguese data. While for the Portugal case study light truncation $H_0^{(2)}$ or Pareto-type behaviour $H_0^{(1)}$ cannot be rejected, the rough truncation model fits the Alberta wildfires data much better than the unbounded model.

   \noindent
  For the Alberta data set, in Figure \ref{fig:2a} (\emph{top right}), the TPa QQ-plot in \eqref{TPQQ}, associated with the validity of \eqref{QT}, has been built on the chosen value $k^*=2232$, which maximizes the correlation between $\log X_{n-j+1,n} $ and $\log\left( \hat D_{T} + j/n\right)$, $j=1,\cdots,k$, for $10<k<n$. For the Pareto index $\alpha$, high quantile  and right endpoint estimation in Figure \ref{fig:2a} we get  conclusions similar to the ones of the earthquake fatalities data set.

\vspace{0.3cm}
The finite sample behaviour of the proposed estimators  $\hat{\alpha}^T_{k,n}$ based on (\ref{Hilltruncate}) and (\ref{onestep}),  $\hat{q}^T_{p,k,n}$ from (\ref{Qest1}), and   $\hat{T}_{k,n}$ from  \eqref{QTinf} has been studied through an extensive Monte Carlo simulation procedure with 1000 runs, both for truncated and non-truncated Pa-type distributions. Here we will only present results concerning Pa and Burr distributions, with truncated and non-truncated versions with sample size $n=400$:
\begin{enumerate}
  \item \emph{Non-truncated models}
  \begin{enumerate}
    \item \emph{Pareto}($\alpha$), $\alpha= 0.5,\; 2$
    \begin{equation}\label{pareto}
F_T(x) = F_W(x)=1-x^{-\alpha}, \,\,\, x>1,\,\,\,\alpha>0,
\end{equation}
    \item \emph{Burr}($\alpha,\rho$),  $\alpha=2$, $\rho=-1$
    \begin{equation}\label{burr}
F_T(x) = F_W(x)=1-(1+x^{-\rho \alpha})^{1/\rho}, \,\,\, x>0,\,\, \rho<0, \,\,\,\alpha>0.
\end{equation}
  \end{enumerate}

  \item \emph{Truncated models}
    \begin{enumerate}
    \item \emph{Truncated-Pareto}($\alpha,T$), $\alpha=0.5,\; 2$ and $T$ a high quantile from the corresponding Pareto model \eqref{pareto}
    \begin{equation}\label{trunc-pareto}
F_T(x)=\frac{1-x^{-\alpha}}{1-T^{-\alpha}}, \,\,\, 1<x<T,\,\,\,\alpha>0.
\end{equation}
    \item \emph{Truncated-Burr}($\alpha,\rho,T$),  $\alpha=2$, $\rho=-1$ and $T$ a high quantile from the corresponding Burr model in \eqref{burr}
    \begin{equation}\label{trunc-burr}
F_T(x)=\frac{1-(1+x^{-\rho \alpha})^{1/\rho}}{1-(1+T^{-\rho \alpha})^{1/\rho}}, \,\,\, 0<x<T,\,\, \rho<0, \,\,\,\alpha>0.
\end{equation}
Note that in case of \eqref{burr} and (\ref{trunc-burr}), $\ell_U (y) = 1+ y^{\rho}(\alpha \rho)^{-1}(1+o(1))$ when $y \to \infty$.
  \end{enumerate}
\end{enumerate}

\noindent
For a particular data set from an unknown but apparently heavy-tailed distribution, the  practitioner does not know if the distribution   comes from a truncated or a non-truncated Pa-type distribution and hence we have to study the behaviour of the tests and the proposed estimators under both cases, and compare them with the existing estimators where appropriate. Our simulation results illustrate this, using three columns in Figures 5-7, setting the cases where $T$ equals the  90 and 99 percentile of the corresponding non-truncated Pa random variable $W$ next to the case where $T=\infty$. The case $T= Q_W(0.90)$ serves as an illustration of the rough truncation case, while $T= Q_W(0.99)$ is meant to represent light truncation.

\begin{figure}[H]
  \centering
 \includegraphics[width=0.32\textwidth]{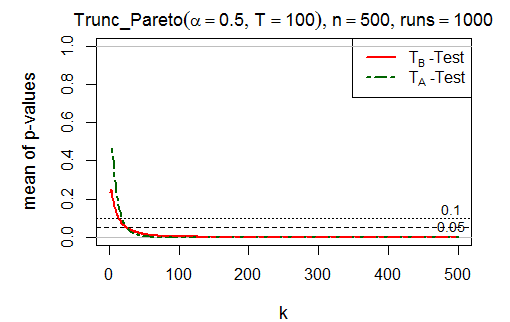} \includegraphics[width=0.32\textwidth]{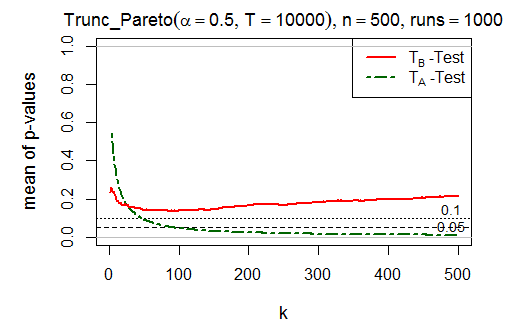} \includegraphics[width=0.32\textwidth]{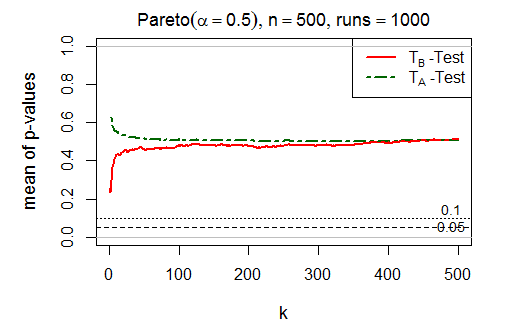}\\
  \includegraphics[width=0.32\textwidth]{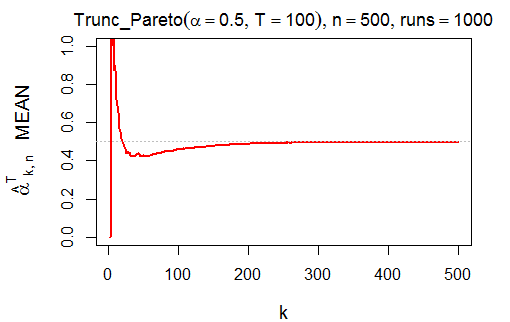} \includegraphics[width=0.32\textwidth]{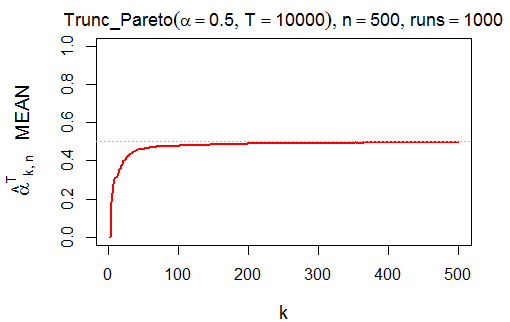} \includegraphics[width=0.32\textwidth]{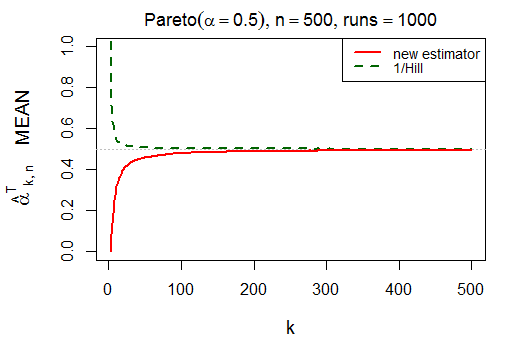}\\
   \includegraphics[width=0.32\textwidth]{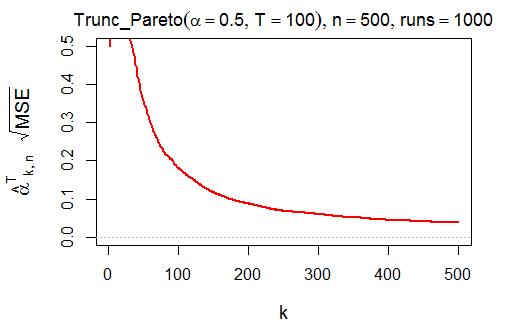} \includegraphics[width=0.32\textwidth]{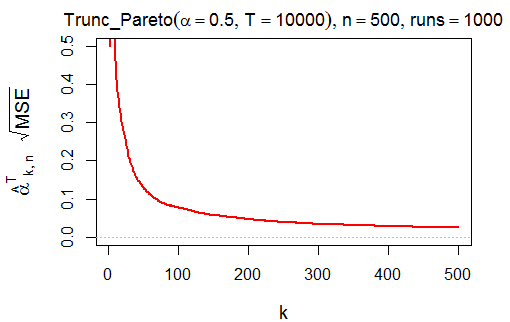} \includegraphics[width=0.32\textwidth]{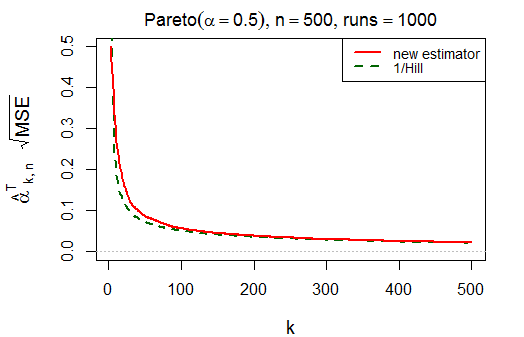}\\
   \includegraphics[width=0.32\textwidth]{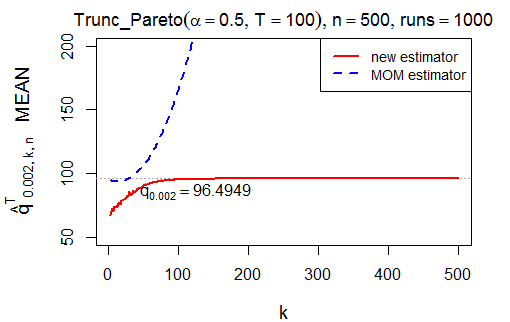} \includegraphics[width=0.32\textwidth]{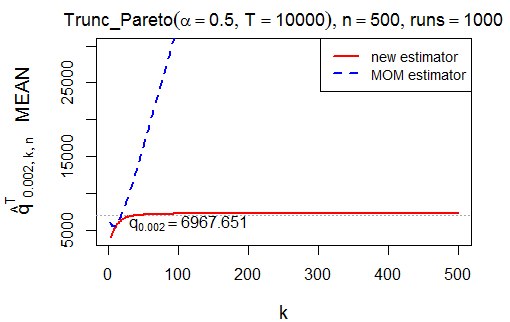} \includegraphics[width=0.32\textwidth]{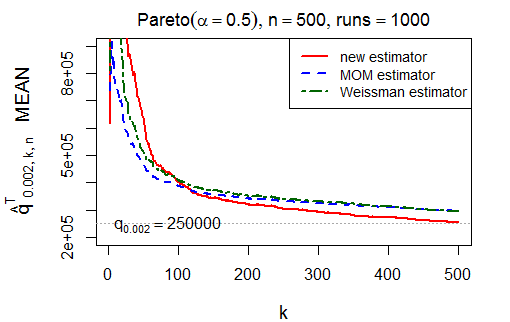}\\   \includegraphics[width=0.32\textwidth]{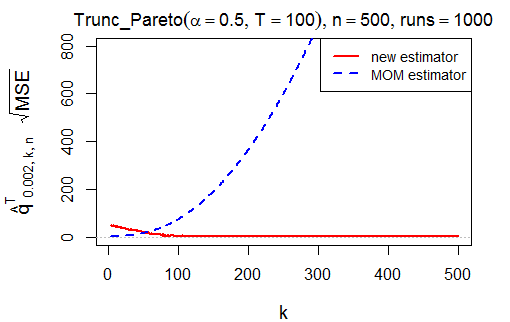} \includegraphics[width=0.32\textwidth]{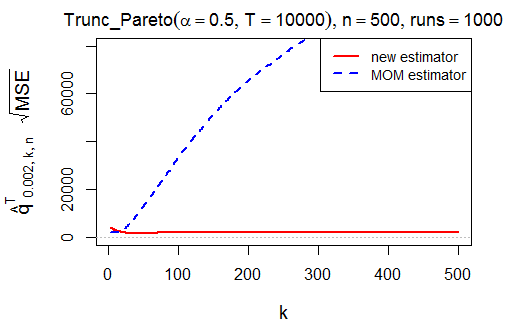} \includegraphics[width=0.32\textwidth]{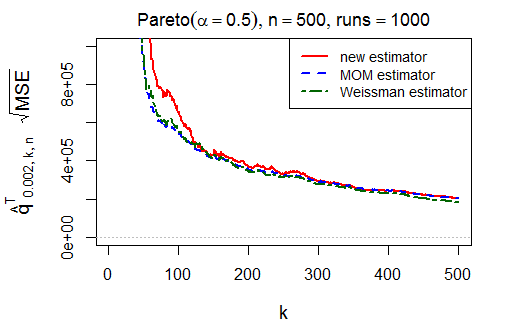}
  \caption{\small\scriptsize Pa$(\alpha=0.5)$: Left column: $T=Q_W(0.90)$; Middle column:  $T=Q_W(0.99)$; Right column:  $T=\infty$. Means of p-values of tests based on $T_A$ and $T_B$ (first row). Estimation of $\alpha$  using the  Newton-Raphson procedure with initial value $\hat\alpha^{(0)} = 1/ H_{k,n}$: mean (second row) and $\sqrt{\mbox{MSE}}$ (third row). Estimation of the high quantile $q_{0.002}$  using  $\hat q^T_{0.002,k,n}$,  $\hat q^{MOM}_{0.002,k,n}$, and $\hat q^W_{0.002,k,n}$ (last column): means (fourth row) and $\sqrt{\mbox{MSE}}$ (fifth row).}
\label{Pareto1}
\end{figure}

\noindent
The mean of the p-values show that both tests $T_A$ and $T_B$ strongly reject the hypotheses $H_0^{(1)}$ and $H_0^{(2)}$ in case
$T=Q_W (0.90)$. In the case $T=Q_W (0.99)$ test $T_A$ rejects more clearly in contrast to the test $T_B$ which then appears more appropriate to test $H_0^{(2)}$. In case of a Burr distribution for $W$ with $T=Q_W (0.99)$, $T_B$ rejects $H_0^{(2)}$ more often for $k \geq 100$ as the deviation from the Pa-model begins to set in for this Pa-type distribution $W$.

\begin{figure}[!ht]
  \centering
 \includegraphics[width=0.32\textwidth]{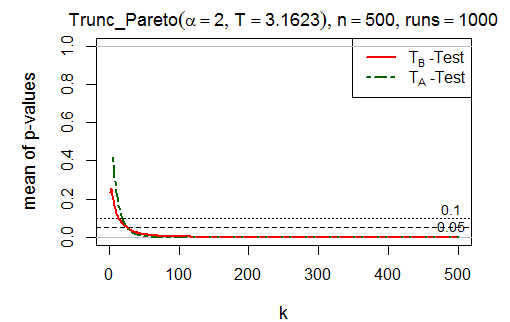} \includegraphics[width=0.32\textwidth]{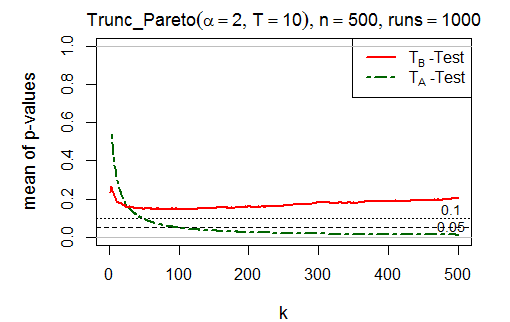} \includegraphics[width=0.32\textwidth]{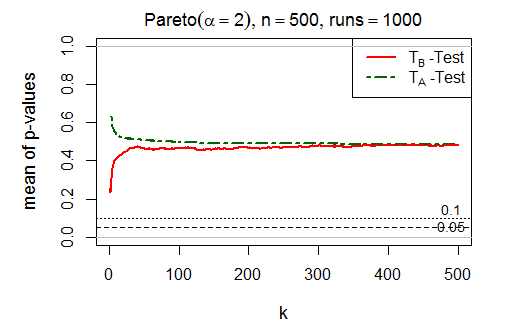}\\
  \includegraphics[width=0.32\textwidth]{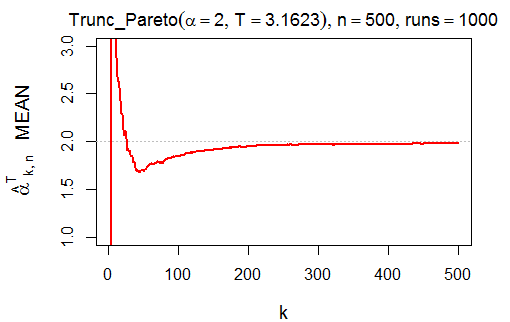} \includegraphics[width=0.32\textwidth]{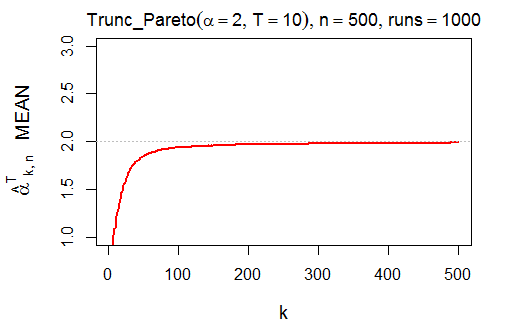} \includegraphics[width=0.32\textwidth]{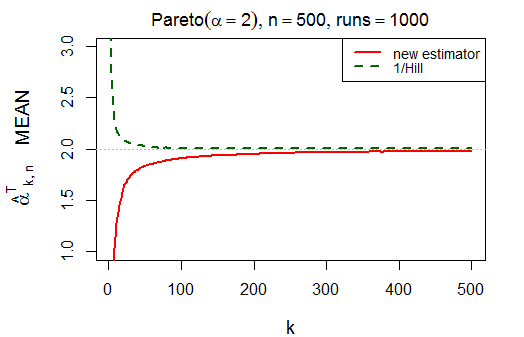}\\
   \includegraphics[width=0.32\textwidth]{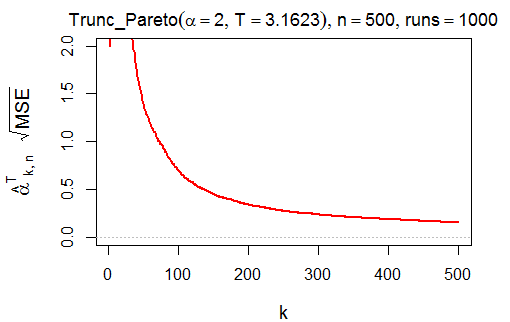} \includegraphics[width=0.32\textwidth]{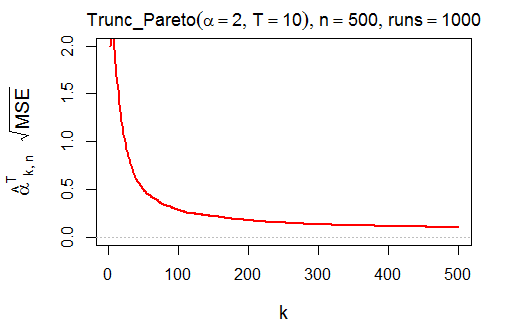} \includegraphics[width=0.32\textwidth]{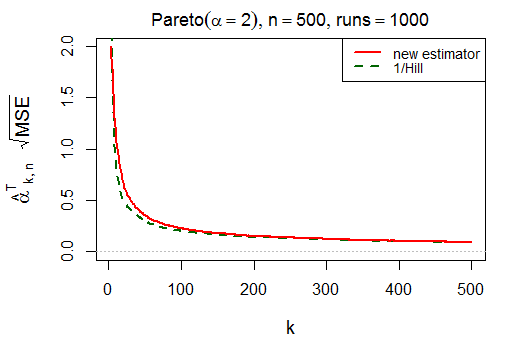}\\
   \includegraphics[width=0.32\textwidth]{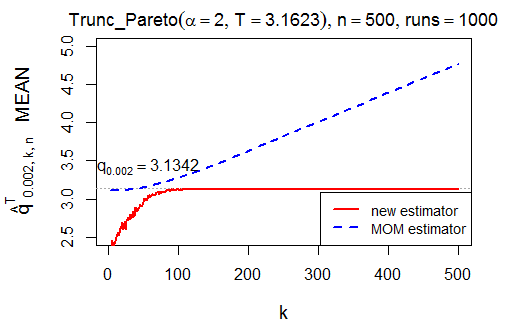} \includegraphics[width=0.32\textwidth]{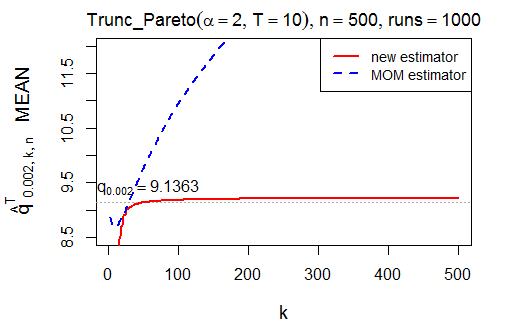} \includegraphics[width=0.32\textwidth]{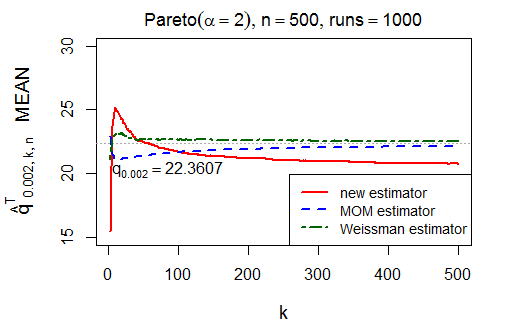}\\   \includegraphics[width=0.32\textwidth]{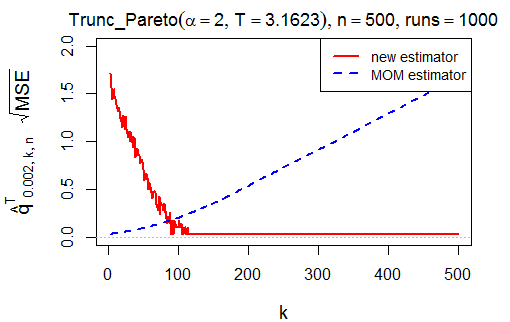} \includegraphics[width=0.32\textwidth]{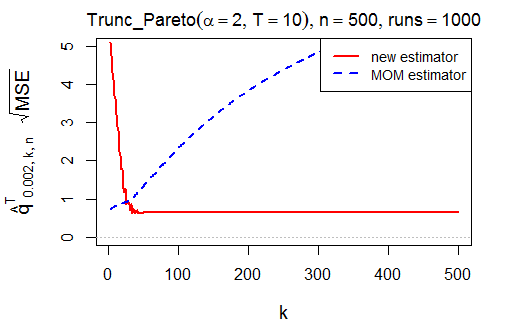} \includegraphics[width=0.32\textwidth]{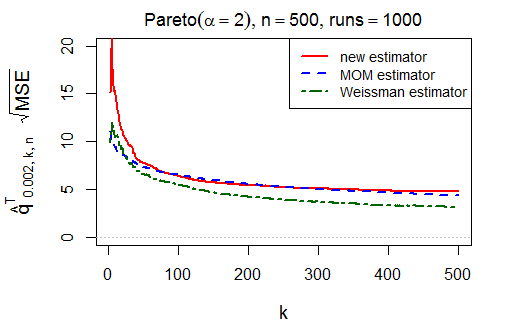}
  \caption{\small\scriptsize Pa$(\alpha=2)$: Left column: $T=Q_W(0.90)$; Middle column:  $T=Q_W(0.99)$; Right column:  $T=\infty$. Means of p-values of tests based on $T_A$ and $T_B$ (first row). Estimation of $\alpha$  using the  Newton-Raphson procedure with initial value $\hat\alpha^{(0)} = 1/ H_{k,n}$: mean (second row) and $\sqrt{\mbox{MSE}}$ (third row). Estimation of the high quantile $q_{0.002}$  using  $\hat q^T_{0.002,k,n}$,  $\hat q^{MOM}_{0.002,k,n}$, and $\hat q^W_{0.002,k,n}$ (last column): means (fourth row) and $\sqrt{\mbox{MSE}}$ (fifth row).}
\label{Pareto2}
\end{figure}
\vspace{0.2cm}\noindent
The  $\alpha$-estimator $\hat{\alpha}$ represented in Figures 5-7 is the solution of \eqref{Hilltruncate}, approximated using the Newton-Raphson iteration as in \eqref{onestep}, with an initial value $\hat\alpha^{(0)} = 1/ H_{k,n}$. With finite samples and fixed $T$, TPa-type distributions belong to the Weibull domain of attraction for maxima with EVI
$\xi=-1$ so that the moment estimator in (\ref{MOM}) almost surely converges to -1. Also, for these models, $1/H_{k,n}$ does not constitute  a consistent estimator either for $\alpha$ or for $\xi$, since in case $\xi <0$ the Hill estimator $H_{k,n}$ almost surely tends to zero when $k/n\rightarrow 0$ as $k,n\rightarrow \infty$. Only when $T=\infty$ we have that $\hat{\alpha}^T_{k,n}$ and $1/H_{k,n}$ estimate the same value $1/\xi$. In case of the truncated Burr distribution the  Pareto
index estimator $\hat{\alpha}^T_{k,n}$  is underestimating $\alpha$, which is not uncommon in extreme value analysis, and is in fact comparable to the behaviour of $1/H_{k,n}$ in this case.

\vspace{0.2cm}\noindent
When estimating an extreme quantile the  estimator in (\ref{quantile-MOM}) based on the moment estimator is designed both for truncated and non-truncated cases and is to be compared with the estimation procedure defined in (\ref{QTinf}). Finally $\hat{q}^T_{p,k,n}$ and  the Weissman (1978) extreme quantile estimator $\hat{q}_{p,k,n}^{W}$ from (\ref{W}) are competitors in case of non-truncated Pa-type distributions only.\\
The convergence of the new quantile estimators seems to be attained at  low thresholds (or high $k$) with high accuracy, contrasting with higher thresholds (or low $k$)  for MOM class estimators. With quantile estimation   an erratic behaviour appears under rough truncation in Figures 6-7 for some smaller  values of $k$. This is a consequence of the use of $\hat{D}_T^{(0)}=\max \{\hat{D}_T,0\}$ rather than $\hat{D}_T$ in practice. If we assume that $T$ is finite then using simply $\hat{D}_T$ rather than $\hat{D}_T^{(0)}$ produces much smoother performance in extreme quantile estimation.
On the other hand in case of non-truncated models the use of $\hat{D}_T$ instead of  $\hat{D}_T^{(0)}$, leads to extreme quantile estimates that are quite sensitive with respect to the value of $D_T$. While the stable parts in the plots of quantile estimates are readily apparent anyway,  we here use    $\hat{D}_T^{(0)}$ in (\ref{Qest1}).
\\\\
In case of non-truncated Pa-type models (right columns) concerning high quantile estimation,
taking into account that $\hat q_p^{MOM}$ and $\hat q_p^{W}$ are designed for this particular situation, we can conclude  that the newly proposed estimators  perform reasonably well at $p$ around $1/n$ if we compare with the classical Weissman and moment-type extreme quantile estimators. For instance in case of the Burr distribution in Figure 7 it appears that our quantile estimator is slightly worse than the moment estimator but better than the Weissman (1978) estimator.  Finally note that for quantile estimators the relative $\sqrt{MSE}$ values can be obtained dividing the presented absolute error $\sqrt{MSE}$ by the exact value of $Q(1-p)$.

\begin{figure}[!ht]
  \centering
 \includegraphics[width=0.32\textwidth]{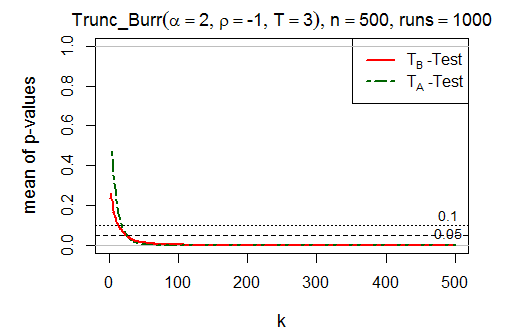} \includegraphics[width=0.32\textwidth]{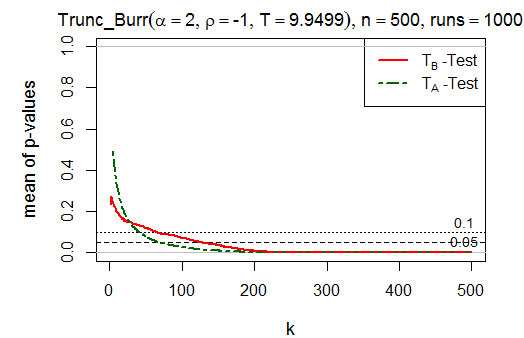} \includegraphics[width=0.32\textwidth]{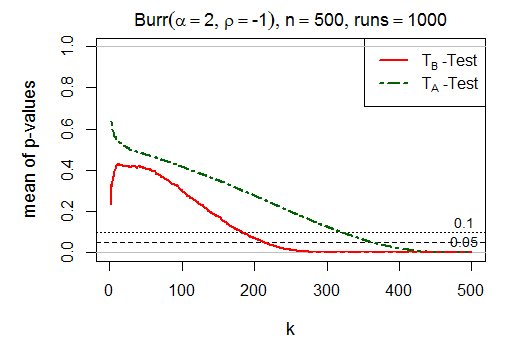}\\
  \includegraphics[width=0.32\textwidth]{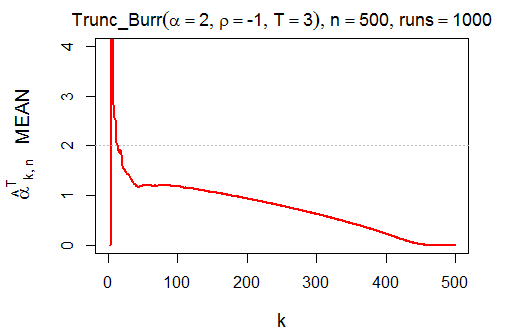} \includegraphics[width=0.32\textwidth]{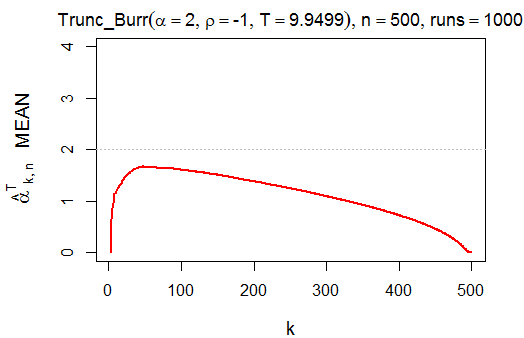} \includegraphics[width=0.32\textwidth]{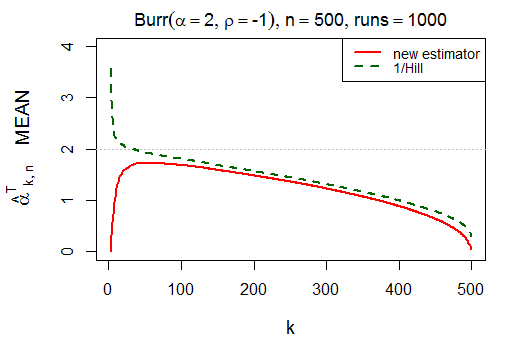}\\
   \includegraphics[width=0.32\textwidth]{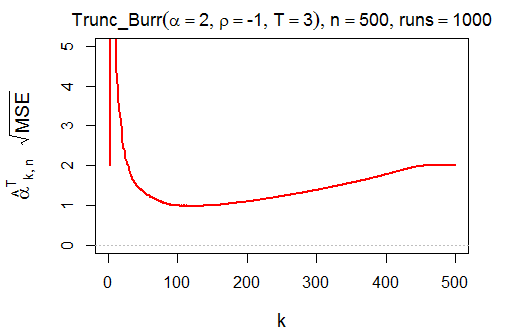} \includegraphics[width=0.32\textwidth]{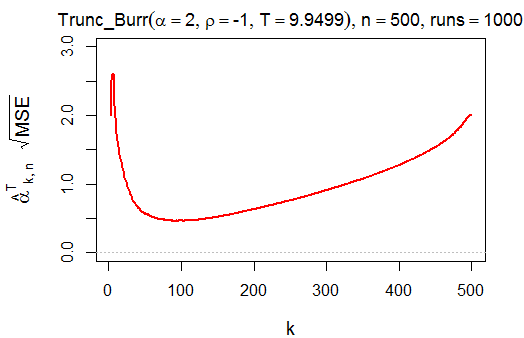} \includegraphics[width=0.32\textwidth]{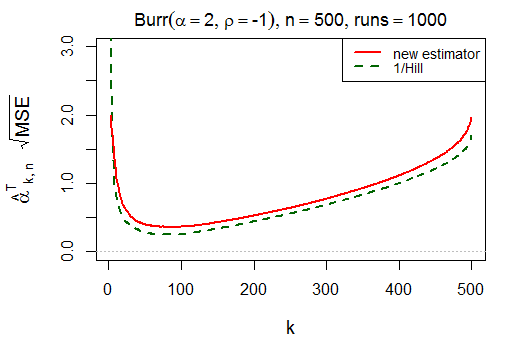}\\
   \includegraphics[width=0.32\textwidth]{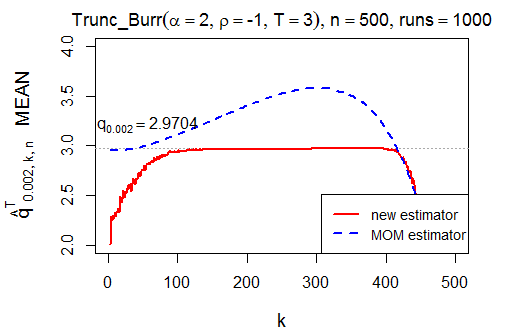} \includegraphics[width=0.32\textwidth]{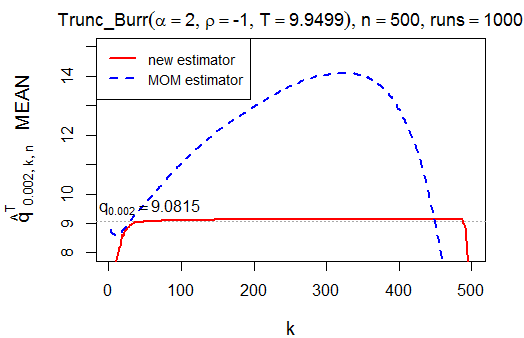} \includegraphics[width=0.32\textwidth]{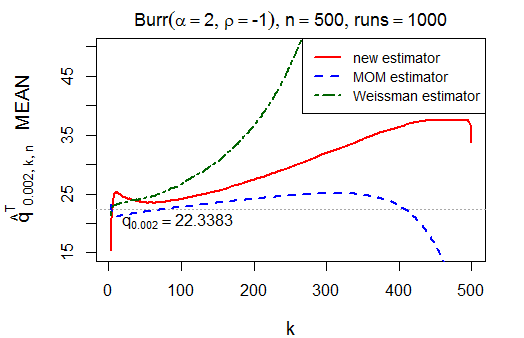}\\   \includegraphics[width=0.32\textwidth]{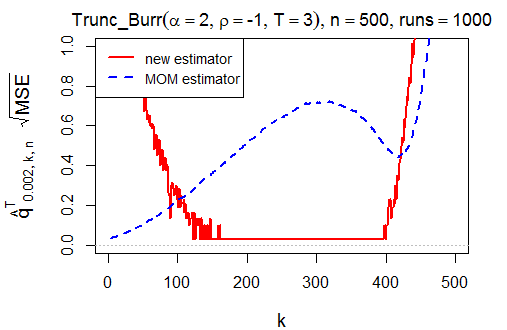} \includegraphics[width=0.32\textwidth]{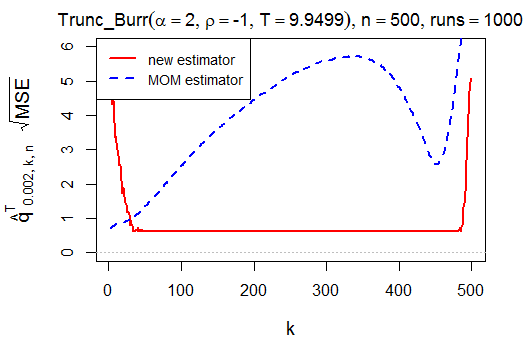} \includegraphics[width=0.32\textwidth]{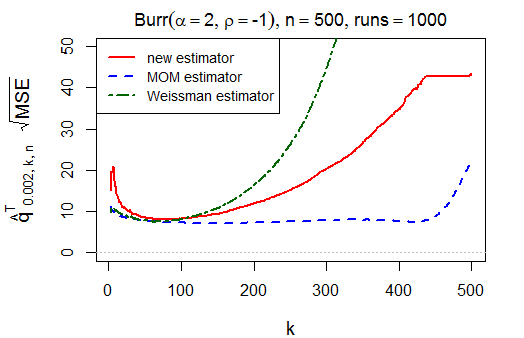}
  \caption{\small\scriptsize Burr$(\alpha=2, \rho = -1)$: Left column: $T=Q_W(0.90)$; Middle column:  $T=Q_W(0.99)$; Right column:  $T=\infty$. Means of p-values of tests based on $T_A$ and $T_B$ (first row). Estimation of $\alpha$  using the  Newton-Raphson procedure with initial value $\hat\alpha^{(0)} = 1/ H_{k,n}$: mean (second row) and $\sqrt{\mbox{MSE}}$ (third row). Estimation of the high quantile $q_{0.002}$  using  $\hat q^T_{0.002,k,n}$,  $\hat q^{MOM}_{0.002,k,n}$, and $\hat q^W_{0.002,k,n}$ (last column): means (fourth row) and $\sqrt{\mbox{MSE}}$ (fifth row).}
\label{Burr}
\end{figure}

\section{Asymptotics of  estimators and tests}

In this section we state the large sample distribution of $\hat{\alpha}^T_{k,n}$ and $\hat{q}^T_{p,k,n}$ defined in (\ref{Hilltruncate}) and \eqref{Qest1}   for TPa-type distributions both under rough \eqref{FT} and light \eqref{LT} truncation.
To this end we will make use of the expressions \eqref{QT}-\eqref{zeta} and \eqref{QTb} for the upper quantile function $Q_T(1-p)$ both under rough and light truncation assuming that $F_T$ is continuous. We also make use of a second order slow variation condition on
 $\ell_U$ specifying the rate of convergence of $\ell_U(tx)/\ell_U(x)$ to 1 as $x\to \infty$, which is used typically in all asymptotic results in extreme value methods (see for instance Theorem 3.2.5 in de Haan and Ferreira, 2006):
\begin{equation}
\lim_{x \to \infty}{1 \over b(x)}\log \frac{\ell_U(tx)}{\ell_U(x)} =  h_{\rho} (t)
\label{Hall3}
\end{equation}
with $\rho <0$, $h_{\rho}(t) = (t^{\rho}-1)/\rho$, and $b$ regularly varying with index $\rho$, i.e. $b(tx)/b(x) \to t^{\rho}$ as $x \to \infty$ for every $t >0$.
Finally $\mathcal{W}$ represents a standard Wiener process.

\begin{thm}\label{theo1} Let \eqref{Hall3} hold and let $n, k=k_n \to \infty$, $k/n \to 0$, $T \to \infty$. Then
\begin{enumerate}[(a)]
  \item if $k/(nD_T) \to \kappa \in (0, \infty )$,
\begin{equation*}
1/\hat{\alpha}^T_{k,n} - 1/\alpha =  {1 \over \delta_{\kappa}}  \left( {1 \over \alpha \sqrt{k}}  \mathcal{N}_{\kappa}^{(1)}  +
b(1/\bar{F}_W(T))\beta_{\kappa} \right)(1 +o_p(1)),
\end{equation*}
with asymptotic variance $1/(k\delta_{\kappa}\alpha^2)$ and
\begin{eqnarray*}
\mathcal{N}_{\kappa}^{(1)} &=& \mathcal{W}(1)  \left\{ 1- {1 \over \kappa }\log \left(1+\kappa \right)  \right\} -  \int_{0}^1
 \mathcal{W}(u) d\log \left( 1+\kappa u\right), \\
\beta_{\kappa} &=& A_{\kappa} - B_{\kappa}c_{\kappa},\\
A_{\kappa}&=& \int_{0}^{1} h_{\rho} ([1+\kappa u]^{-1} )du -  h_{\rho} ([1+\kappa ]^{-1} ),\\
B_{\kappa} &=& h_{\rho} ([1+\kappa]^{-1} ),\\
c_{\kappa} &=& {1 \over \kappa } - { (1+ \kappa ) \over \kappa^2 }
\log \left(  1+ \kappa \right), \\
\delta_{\kappa}&=& 1- {1+ \kappa \over \kappa^2 }
 \log^2\left( 1+ \kappa  \right);
\end{eqnarray*}
  \item if $k/(nD_T) \to \infty$ and $D_T \log (k/(nD_T)) = o((n/k)^{-1+ \rho})$
\begin{equation*}
1/\hat{\alpha}^T_{k,n} - 1/\alpha =  \left( {1 \over \alpha \sqrt{k}} \mathcal{N}^{(2)} + b(n/k) (1 - \rho)^{-1}\right)(1+ o_p(1)),
\end{equation*}
with
$\mathcal{N}^{(2)} = \mathcal{W}(1) -\int_{0}^1 (\mathcal{W}(u)/u) \; du \sim \mathcal{N}(0,1)$.
\end{enumerate}
\end{thm}

\vspace{0.2cm}\noindent
{\bf Remark 1.}  In case $k/(nD_T) \to \infty$ the asymptotic result for $1/\hat{\alpha}^{T}_{k,n}$ is identical to that of the Hill estimator $H_{k,n}$ under a Pa-type distribution, as given for instance in Beirlant {\it et al.} (2004), section 4.2.

\begin{thm} \label{theo2}  Suppose \eqref{Hall3} holds and  $n, k=k_n \to \infty$, $k/n \to 0$, $T \to \infty$  and $p=p_n$ such that $np_n = o(k)$. Moreover $E$ denotes a standard exponential random variable.
\begin{enumerate}[(a)]
  \item Let $k/(nD_T) \to \kappa$. Then
\begin{eqnarray*}
\log \hat{q}^T_{p,k,n} - \log q_p &=&  \left\{
-\kappa \left({p(n+1)\over k+1}-{1 \over k+1}\right)\left( {1 \over \alpha^T_{k,n}}- {1 \over \alpha}\right)\right. \\
&& \hspace{0.5cm} - {\kappa \over \alpha (k+1)}(E -1) \\
&& \hspace{0.5cm} \left.
+\kappa b(1/\bar{F}_W(T))\left({p(n+1)\over k+1}-{1 \over k+1} \right)\right\}
(1+o_p(1)).
\end{eqnarray*}
  \item Let $np_n \to \infty$, $\log(np_n) = o(\sqrt{k})$, $\sqrt{k} nD_T/ \log (k/np_n) \to 0$, and $b(n/k) \log (k/(np_n) \to 0$. Then
\begin{eqnarray*}
&& \hspace{-1cm} \log \hat{q}_{p,k,n} - \log q_p  \\
&=& \log \left( k/(np_n) \right)
\left\{ {1 \over \alpha \sqrt{k}}\mathcal{N}^{(2)}
+ b(n/k) (1-\rho)^{-1}
\right\}(1+ o_p(1)) \\
&& -{1 \over \alpha}  {1 \over np_n}(E-1+o_p(1)).
\end{eqnarray*}
If $np_n \to c>0$ then
$$
\log \hat{q}_{p,k,n} - \log q_p =  -{1 \over \alpha}  {1 \over c}(E-1)(1+o_p(1)).
$$
\end{enumerate}
\end{thm}

 \noindent
{\bf Remark 2.} In case $k/(nD_T) \to \kappa$ both the asymptotic bias and the stochastic part of $\hat{q}^T_{p,k,n}$ are of smaller order than in case of $1/\hat{\alpha}^T_{k,n}$. This  confirms the findings of the simulations where the
plots of the quantile estimators are found to be quite horizontal as a function of $k$ and show a small variance, compared to other extreme quantile estimators in this case.

\vspace{0.2cm} \noindent
In case $k/(nD_T) \to \infty$, note that the quantile estimator $\hat{q}^T_{p,k,n}$ is only consistent if $(np_n)^{-1} \to 0$, this is for quantiles $q_p$  situated maximally up to the border of the sample $q_{n^{-1}}$, using for instance a sequence of the type  $p_n = (\log k)^\tau /n$ for some $\tau >0$. The extra factor
$\left( {(1+ { n\hat{D}_{T} \over k}) /( 1+ { \hat{D}_{T}\over p})} \right)^{1/\hat{\alpha}^T_{k,n}}$ in (\ref{QTinf}) compared to the Weissman estimator $\hat{q}^W_{p,k,n}$ induces this restriction. The first term in the expansion of $\log\hat{q}^T_{p,k,n}$ in Theorem 2(b) is indeed the asymptotic expansion of  $\hat{q}^W_{p,k,n}$ as given for instance in Beirlant {\it et al.} (2004), section 4.6. If $np_n\to \infty$ and ${\sqrt{k} \over np_n \log((k+1)/((n+1)p_n)) } \to 0 $ then the expansion of the Weissman estimator is dominant, while if $np_n\to \infty$ and ${\sqrt{k} \over np_n \log(k/(np_n)) } \to \infty $ the second term in the expansion is to be retained.
\\\\
From Theorem 2 it follows that it is important to be able to test if for a given case study rough or light truncation holds. Indeed, if rough truncation holds then in extreme quantile estimation  $\hat{q}^T_{p,k,n}$ from \eqref{Qest1} should be used, while under light truncation the estimator $\hat{q}^{\infty}_{p,k,n}$ from \eqref{QTWT}, or  the classical Weissman estimator \eqref{W}, or the moment type extreme quantile estimator $\hat{q}^{MOM}_p$ can be used when extrapolating outside the sample. From Theorem 3(a) the consistency of both tests \eqref{Atest} and \eqref{Btest} follows directly, while the null distributions are conisdered in parts (b) and (c).
\begin{thm} \label{theo3}  Suppose \eqref{Hall3} holds and  $n, k=k_n \to \infty$, $k/n \to 0$, $T \to \infty$.
\begin{enumerate}[(a)]
  \item Let $k/(nD_T) \to \kappa$. Then
\begin{eqnarray*}
{\log R_{k,n} \over H_{k,n}} &\to_p & {\kappa \log (1+\kappa) \over \log (1+\kappa )-\kappa}, \\
L_{k,n}(H^{-1}_{k,n}) &\to_p& { \{ \kappa - \log (1+\kappa)\}(1+\kappa)^{-{\kappa \over \kappa - \log (1+\kappa)}}
\{ (1+\kappa)^{{2\kappa - \log (1+\kappa) \over \kappa - \log (1+\kappa) }} -1 \} \over \kappa (2\kappa - \log (1+\kappa) )} .
\end{eqnarray*}
  \item Let  $nD_T \to 0$. Then
\begin{equation*}
T_{A,k,n} \to_d E
\end{equation*}
where $E$ is a standard exponential random variable.
\item Let $nD_T/k \to 0$ and $\sqrt{k}(n/k)^{\rho} = o(1)$. Then
$$
T_{B,k,n} \to_d \mathcal{N}(0,1).
$$
\end{enumerate}
\end{thm}
\noindent
Note that the result in Theorem 3(b) needs a stronger condition on $T$ for the limit $E$ to hold under $H_0^{(2)}$, meaning that for the test based on $T_{A}$ the truncation point $T$ has to lie higher than in case of $T_B$ in order to keep a given significance level. This yields a theoretical confirmation of the simulation results where the $T_A$ test was found to reject light truncation situations that deviate from the untruncated Pa-type distributions, sooner than the $T_B$ test.

\section{Conclusion}

We have extended the work on estimating the Pareto index $\alpha$ under truncation from Aban {\em et al.} (2006) and Nuyts (2010) to extreme quantile estimation, and considered also truncation of regularly varying tails.  The main proposals and findings are
\begin{itemize}
\item The new estimator of the Pareto index $\alpha$ is effective whether the underlying distribution is truncated or not,  thus unifying previous approaches. Although based on a truncated model, the  estimator of $\alpha$ is competitive even when the underlying distribution is unbounded.
\item Our method leads to new quantile estimators which are especially effective in the case of rough truncation.  In case the data come from a light truncated Pa-type distribution, which is the case when the Pa QQ-plot \eqref{PaQQ} is linear in the right tail, the extreme quantile estimator \eqref{Qest1} should not be used for extrapolation far out of range of the available observations as discussed in Remark 2.
\item  A new TPa QQ-plot is constructed that can assist in verifying the validity of the TPa-type model. Moreover a new test is provided for testing light truncation against rough truncation which offers an extra
tool for practitioners.
\end{itemize}

\noindent
\textbf{Acknowledgment} The authors thank Mark Meerschaert for helpful discussions and many suggestions that had a significant positive influence on the paper.

\pagebreak

\begin{center}
\textbf{\LARGE \textmd{Tail fitting for truncated and non-truncated
Pareto-type distributions}
\\\vspace{0.5cm}
{\it\small Supplementary material: Proofs of the asymptotic results}}\\\vspace{0.5cm}
\author{ Beirlant J.$^{a}$\footnote{Corresponding author: Jan Beirlant, KU Leuven, Dept of Mathematics and  LStat, Celestijnenlaan 200B, 3001 Heverlee, Belgium; Email: jan.beirlant@wis.kuleuven.be },  Fraga Alves, M.I.$^{b}$,  Gomes, M.I.$^{b}$, \\
{$^a$ \fontsize{8pt}{11pt} \selectfont Department of Mathematics and Leuven Statistics Research Center, KU Leuven}
\\
{$^b$ \fontsize{8pt}{11pt} \selectfont Department of Statistics and Operations Research, University of Lisbon }
 }
\end{center}
\maketitle
\setcounter{prop}{0}
\setcounter{equation}{0}
\setcounter{page}{1}

\makeatletter
\renewcommand{\theequation}{S\arabic{equation}}
\makeatother
\makeatletter
\renewcommand*{\@biblabel}[1]{[S#1]}
\makeatother
\makeatletter
\renewcommand*{\@cite}[1]{[S#1]}
\makeatother

\noindent
{\small
\textbf{Proof of Theorem 1}
The mean value theorem implies that $1/\hat{\alpha}^T_{k,n} - 1/\alpha = -f(1/\alpha)/f'(1/\tilde{\alpha})$ where $\tilde{\alpha}=\tilde{\alpha}^T_{k,n}$ is between $\alpha$ and $\hat{\alpha}^T_{k,n}$, with $f({1 \over \alpha})=H_{k,n} - 1/\alpha - {R_{k,n}^{\alpha} \over 1-R_{k,n}^{\alpha} }\log R_{k,n}$.
Then  that the limit distribution of $1/\hat{\alpha}^T_{k,n}$ is found from the asymptotic distribution of
$$
\left( 1- \tilde{\alpha}^{2} {R_{k,n}^{\tilde\alpha}\log^2 R_{k,n} \over (1- R_{k,n}^{\tilde\alpha})^2}  \right)^{-1}
\left( H_{k,n}-1/\alpha - {R_{k,n}^{\alpha}\log R_{k,n} \over 1- R_{k,n}^{\alpha}} \right).
$$
Hence the asymptotic behaviour of $H_{k,n}$ and $\log R_{k,n}$ constitute essential building blocks in the derivation of the asymptotics for $\hat{\alpha}^T_{k,n}$. We consider these in the following Propositions.

\vspace{0.2cm} \noindent
For this we  make use of the result (see de Haan and Ferreira, 2006, 7.2.12) that for some standard Wiener process $\mathcal{W}$ (with $E(\mathcal{W}(s)\mathcal{W}(t)) = \min (s,t)$) we have uniformly over all $j=1,\ldots,k$, as $k,n \to \infty$, $k/n \to 0$
\begin{equation}
\sqrt{k} \left( {n \over k} U_{j,n} - {j \over k} \right) - \mathcal{W}\left( {j \over k}\right) \to_p 0.
\label{wiener}
\end{equation}

\begin{prop}\label{prop2}
 Let (39) hold and let $n, k=k_n \to \infty$, $k/n \to 0$, $T \to \infty$. Then
 \begin{enumerate}[(a)]
   \item if $k/(nD_T)\to \kappa$,
\begin{eqnarray*}
H_{k,n} && \\ &&  \hspace{-1.5cm} = {1 \over \alpha}
\left(1- {1 \over k/(nD_T) }  \log \left( 1+ {k \over nD_T} \right) \right) \\
 && \hspace{-1cm}
 + {1 \over \alpha \sqrt{k}} \left( {\mathcal{W}(1) \, k/(nD_T) \over 1+k/(nD_T)}-  \int_{0}^{1} \mathcal{W}(u) d\log (1+{k \over nD_T}u) \right)(1+o_p(1)) \\
 && \hspace{-1cm} +  b(1/\bar{F}_W(T)) A_{k,n,T}(1+o_p(1)),
\end{eqnarray*}
where
$$A_{k,n,T} =
 \int_{0}^{1} h_{\rho} ([1+{k \over nD_T} u]^{-1} )du -  h_{\rho} ([1+{ k \over nD_T}  ]^{-1} ),
$$
   \item  if $k/(nD_T) \to \infty$,
\begin{eqnarray*}
H_{k,n} & = & {1 \over \alpha} + {1 \over \alpha \sqrt{k}} \left( \mathcal{W}(1) - \int_{0}^{1} {\mathcal{W}(u) \over u } du \right) (1+o_p(1)) \\
 && + b(n/k)
   \int_{0}^{1}
  h_{\rho}(u^{-1})du \, (1+o_p(1)) \\
&& +  {1 \over \alpha}{nD_T\over k}\log ({k \over nD_T}) \, (1+o_p(1)) ,
\end{eqnarray*}
where the first two terms in this expansion are the limits for $k/(nD_T) \to \infty$ of the first two lines in the expansion in case  (a).
 \end{enumerate}
\end{prop}
\vspace{0.2cm}\noindent
\begin{proof}
Let $j_r =j-r+1$
and let $U_{1,n} \leq U_{2,n} \leq \ldots \leq U_{n,n}$ denote the  order statistics from an i.i.d. sample of size $n$ from the uniform (0,1) distribution. Then using summation by parts and the fact that $X_{n-j+1,n}=_d Q_T(1-U_{j,n})$ ($j=1,\ldots,n$)
\begin{eqnarray*}
H_{k,n} &=& {1 \over k} \sum_{j=1}^k j \left( \log X_{n-j+1,n}-\log X_{n-j,n}\right) \\
&=& -{1 \over k} \sum_{j=1}^k j \left( \log Q_T (1-U_{j+1,n})-\log Q_T (1-U_{j,n})\right)
\\
&=& - {k+1 \over k}\sum_{j=1}^k {j \over k+1} \int_{U_{j,n}}^{U_{j+1,n}} d \log Q_T (1-w).
\end{eqnarray*}
Using (\ref{wiener}), $H_{k,n}$ can now be approximated as $k,n \to \infty$ by the integral
\begin{equation}
 -k \int_{0}^{1} u
  \left\{ \int_{u+ \mathcal{W}(u)/\sqrt{k}}^{u+\mathcal{W}(u)/\sqrt{k}+ 1/k} d \log Q_T (1-{k \over n}w) \right\} du.
  \nonumber
  \end{equation}
  Using the mean value theorem on the inner integral between $u+ \mathcal{W}(u)/\sqrt{k}$ and $u+ \mathcal{W}(u)/\sqrt{k} + 1/k$, followed by an integration by parts, we obtain the approximation
 \begin{eqnarray}
  & &- \int_{0}^{1} u d \log Q_T \left(1- {k \over n}\left(u + {\mathcal{W}(u)\over \sqrt{k}}\right) \right)  \nonumber \\
 &=&  -\log U_T \left({n \over k}/\left(1 + {\mathcal{W}(1) \over \sqrt{k}}\right)\right)
+\int_{0}^{1}\log U_T \left({n \over k}/\left(u + {\mathcal{W}(u) \over \sqrt{k}}\right)\right) du. \nonumber \\
&&  \label{midH}
  \end{eqnarray}

  \vspace{0.2cm} \noindent
 {\it  First, let $k/(nD_T) \to \kappa$}. Then from (15) the approximation (\ref{midH}) of $H_{k,n}$ equals
  \begin{eqnarray*}
  && {1 \over \alpha} \log \left(1+{k \over nD_T}\left(1 + {\mathcal{W}(1) \over \sqrt{k}}\right)\right) \\
  && - \log \ell_U \left( (1/\bar{F}_W (T)) \left[1+{k \over nD_T}\left(1 + {\mathcal{W}(1) \over \sqrt{k}}\right) \right]^{-1}\right) \\
  && -{1 \over \alpha}
    \int_{0}^{1}\log \left( 1+{k \over nD_T}\left(u + {\mathcal{W}(u) \over \sqrt{k}}\right) \right) du \\
    && + \int_0^{1}\log \ell_U \left( (1/\bar{F}_W (T))
    \left[1+{k \over nD_T}\left(u + {\mathcal{W}(u) \over \sqrt{k}}\right) \right]^{-1}\right) du.
\end{eqnarray*}

Next, add and  subtract $\log \ell_U(1/\bar{F}_W (T))$ from  the second and fourth line respectively, and
use the approximations
$$
\log \left( 1+{k \over nD_T}\left(u + {\mathcal{W}(u) \over \sqrt{k}}\right) \right)
= \log \left(1+{k \over nD_T}u \right) + {\mathcal{W}(u) \over \sqrt{k}} {{k \over nD_T}\over 1+{k \over nD_T}u^*},
$$
with $0 < u \leq 1$ and $u^*$ between $u$ and $u+\mathcal{W}(u)/\sqrt{k}$, and
\begin{eqnarray*}
&& \hspace{-2cm} \log \frac{\ell_U \left((1/\bar{F}_W (T))  [1+{k \over nD_T}(u + {\mathcal{W}(u) \over \sqrt{k}}) ]^{-1}\right)}
{\ell_U(1/\bar{F}_W (T)) }  \\
&=&  b(1/\bar{F}_W (T))  h_{\rho}\left([1+{k \over nD_T}u]^{-1} \right)(1+o_p(1)).
\end{eqnarray*}
Finally, using partial integration we have
$$
  \int_{0}^{1} \log \left( 1+{k \over nD_T}u \right) du  =
-1 + \log (1+{k \over nD_T})   + {nD_T \over k} \log ( 1+ {k  \over nD_T})
$$
from which one obtains the stated result in (a).

\vspace{0.2cm} \noindent
 {\it  Secondly, consider $k/(nD_T) \to \infty$}. Then using (16)
the proof follows the lines of proof of the asymptotic normality of the Hill estimator as for instance given in Mason and Turova (1994), except
for the extra term
$$
e_{k,n}:= -{1 \over \alpha}{1 \over k}\sum_{j=1}^k \log \frac{1 + {D_T \over U_{j,n}}}{1+ {D_T \over U_{k+1,n}}}
$$
appearing from the factor $[1+yD_T]^{-1/\alpha}$ in (16). However, using
$U_{j,n}/(j/(n+1)) \to_p 1$, one derives that
$$
e_{k,n}\sim_p  -{1 \over \alpha}{nD_T \over k}\log \left(1+ {k \over nD_T}\right).
$$
Condition $D_T \log (k/(nD_T)) = o((n/k)^{-1+\rho})$ in Theorem 1(b) entails that the bias term due to the factor $(1+D_Ty)^{-1/\alpha}$ in (16) is negligible with respect to the classical bias $b(n/k)/(1-\rho)$ due to the last factor in (16).
 \end{proof}

\vspace{0.2cm} \noindent
 \begin{prop}\label{prop3}  Let (39) hold and let $n, k=k_n \to \infty$, $k/n \to 0$, $T \to \infty$. Then
 \begin{enumerate}[(a)]
   \item if $k/(nD_T) \to \kappa$,
\begin{eqnarray*}
\log R_{k,n} &=&
- {1 \over \alpha}\log \left( 1+ {k \over nD_T}  \right)  -{1 \over \alpha \sqrt{k}}  \mathcal{W}(1) {k/(nD_T)  \over 1+k/(nD_T)} )(1+o_p(1)) \\
&& + b(1/\bar{F}_W (T)) B_{k,n,T}(1+o(1)),
\end{eqnarray*}
where
$$
B_{k,n,T} =h_{\rho} ([1+(k/(nD_T))]^{-1} ),
$$
   \item    if $k/(nD_T) \to \infty$,
\begin{eqnarray*}
\log R_{k,n} &=&
 -{1 \over \alpha}\log k
 -{1 \over \alpha} \left( E_{k,k} - \log k \right) +{1 \over \alpha}\log \frac{1+ {nD_T \over k} k}{1+ {nD_T \over k}}
 \\
&&  -b(n/k)\left( {nD_T \over k}+ h_{\rho}\left( [k^{-1}+ (nD_T/k)]^{-1}\right)\right)(1+o(1)),
\end{eqnarray*}
where $E_{k,k}$ denotes the maximum of a sample of size $k$ from the standard exponential distribution.
  \end{enumerate}
 \end{prop}
 \vspace{0.2cm}\noindent
\begin{proof} The proof of (a) follows similar lines as the proof of Proposition 1(a). Concerning part (b) remark that using (16)
\begin{eqnarray*}
\log R_{k,n}&=& {1 \over \alpha}\log {U_{1,n}\over U_{k+1,n}}
- {1 \over \alpha}\log \frac{1+ {nD_T \over k}}{1+ {nD_T \over k}k} \\
&&+ \log \frac{\ell_U \left( C_T (n/k) [1+ {nD_T \over k}]^{-1}\right)}
{\ell_U \left( C_T (n/k)\, k[1+ k\,{nD_T\over k}]^{-1}\right)}.
\end{eqnarray*}
The result then follows since $\log {U_{k+1,n}\over U_{1,n}} =_d E_{k,k}$, and by using (39).
\end{proof}

\vspace{0.2cm} \noindent
\begin{proof}\textbf{of Theorem 1} ({\it cont'd}).
First we derive the consistency of $\hat{\alpha}^T_{k,n}$ under the conditions of Theorem 1, so that then $\tilde{\alpha} \to_p \alpha$.
Aban {\it et al.} (2006, see A.4) showed that $\tilde{f}(t):= {1 \over t} + {R_{k,n}^{t} \log R_{k,n} \over 1-  R_{k,n}^{t} }- H_{k,n}$ is a decreasing function in $t \in (0, \infty) $. Moreover
$\lim_{t \to \infty}\tilde{f}(t) = - H_{k,n} <0$ and $\lim_{t \to 0}\tilde{f}(t) = -(\log R_{k,n})/2-H_{k,n}$. Showing that asymptotically under the conditions of the theorem $-(\log R_{k,n})/2-H_{k,n} >0$ using Propositions \ref{prop2} and \ref{prop3} in both cases (a) and (b), we have then that there is a unique solution to the equation $\tilde{f}(t) =0$. Note with Propositions \ref{prop2} and \ref{prop3} that for the true value $\alpha$ we have $\tilde{f}(\alpha)= o_p(1)$, since $H_{k,n}$ and ${1 \over \alpha} + {R_{k,n}^{\alpha} \log R_{k,n} \over 1-  R_{k,n}^{\alpha} }$  asymptotically are equal, namely to
${1 \over \alpha}\left( 1-({k\over nD_T})^{-1}\log  ( 1+{k\over nD_T}) \right)$ in case (a),  and
$\alpha^{-1}$ in case (b). So the true value $\alpha$ asymptotically is a solution from which the consistency follows.

\vspace{0.2cm}
Now using Propositions \ref{prop2} and \ref{prop3} we obtain that
\begin{equation}
 1- \tilde{\alpha}^{2} {R_{k,n}^{\tilde\alpha}\log^2 R_{k,n} \over (1- R_{k,n}^{\tilde\alpha})^2} =  \delta_{k,n,T}(1+o_p(1)),
 \label{denom}
\end{equation}
where
$$
\delta_{k,n,T} = 1- {1+ {k \over nD_T}  \over  ({k \over nD_T})^2 }
\log^2 \left( 1+ {k \over nD_T}  \right).
$$
Next, consider $g(H_{k,n},\log R_{k,n})= H_{k,n}-1/\alpha - {R_{k,n}^{\alpha}\log R_{k,n} \over 1- R_{k,n}^{\alpha}}$ with
$$
g(x,y)= x-{1 \over \alpha} -y{e^{\alpha y} \over 1-e^{\alpha y}}.
$$
The Taylor approximation of $g(H_{k,n},\log R_{k,n})$ around the asymptotic expectated value $E_\infty H_{k,n}$ and $E_\infty \log R_{k,n}$ yields
\begin{eqnarray}
&& \hspace{-1cm} g(E_\infty H_{k,n},E_\infty \log R_{k,n}) + (H_{k,n}-E_\infty H_{k,n}) \nonumber \\ &+ &
(\log R_{k,n}- E_\infty \log R_{k,n}){\partial g \over \partial y}(E_\infty H_{k,n},E_\infty \log R_{k,n} ),
\label{Taylor1}
\end{eqnarray}
with, based on Proposition  \ref{prop3},
\begin{equation}
{\partial g \over \partial y}(E_\infty H_{k,n},E_\infty \log R_{k,n} )
= - c_{k,n,T}(1+o(1)),
\label{partg}
\end{equation}
where
$$
c_{k,n,T}  = ({k \over nD_T})^{-1}\left( 1 - { (1+ {k \over nD_T} ) \over  ({k \over nD_T}) }
\log \left(  1+ {k \over nD_T}  \right) \right).
$$
From Propositions \ref{prop2} and \ref{prop3}, (\ref{Taylor1}), (\ref{partg}), and (\ref{denom}) we find  that the stochastic part in the development of $(\hat{\alpha}^{T}_{k,n})^{-1}-\alpha^{-1}$ is given by
$$
 {1 \over \delta_{k,n,T}\alpha \sqrt{k}}
\left( -  \int_{0}^1
 \mathcal{W}(u) d\log \left( 1+{k \over nD_T}u\right)
 +\mathcal{W}(1)
 \left\{ 1- {\log \left( 1+{k \over nD_T} \right)  \over {k \over nD_T}} \right\} \right).
$$
Developing for  $k/(nD_T) \to \kappa$, respectively  $k/(nD_T) \to \infty$, leads to the  stated asymptotic variances in  cases (a) and (b).
\\\\
From (\ref{Taylor1}), (\ref{partg}), and the asymptotic bias expressions in Propositions \ref{prop2} and \ref{prop3} one finds the asymptotic bias expressions of $(\hat{\alpha}^{T}_{k,n})^{-1}$. For instance in case $k/(nD_T)  \to \kappa$ we find that
\begin{eqnarray}
&& \hspace{-1cm} g(E_\infty H_{k,n},E_\infty \log R_{k,n}) \nonumber \\
&= & b(1/\bar{F}_W(T))
\left( A_{k,n,T}  - B_{k,n,T}c_{k,n,T}\right)(1+o(1)).
\label{Einfinf}
\end{eqnarray}
\end{proof}

\vspace{0.2cm} \noindent
\begin{proof} {\textbf{of Theorem 2.}} First consider the case $k/(nD_T) \to \kappa$. Then from (22)
\begin{eqnarray*}
\log q^T_{p,k,n}&=&
\log X_{n-k,n} + {1 \over \alpha^T_{k,n}}\log \frac{1-{1 \over k+1}}{R_{k,n}^{\hat{\alpha}^T_{k,n}}-{1 \over k+1}} -{1 \over \alpha^T_{k,n}}\log \left( 1+ {p \over \hat{D}_T}\right) \\
&=& \log X_{n,n}-{1 \over \alpha^T_{k,n}}
 \log \frac{1-[(k+1)R_{k,n}^{\hat{\alpha}^T_{k,n}}]^{-1}}{1-{1 \over k+1}} \\
 && -{1 \over \alpha^T_{k,n}}
 \log \left( 1+ {p(n+1) \over k+1}  \frac{k+1}{(n+1)\hat{D}_T}
 \right),
\end{eqnarray*}
while
\begin{eqnarray*}
\log X_{n,n} &=& \log U_T (1/U_{1,n}) \\ &=&
\log T - {1 \over \alpha}log \left(1+ {k+1 \over (n+1)D_T}{E_1 \over \bar{E}_{n+1}}{1 \over k+1} \right)\\
&&+
\log {\ell_U \left( (1/\bar{F}_W (T))[1+{U_{1,n}\over D_T}]^{-1}  \right) \over \ell_U \left( 1/\bar{F}_W (T)  \right)}
\end{eqnarray*}
with $\bar{E}_{n+1}$ the average of an i.i.d. sample $E_1, \ldots, E_n, E_{n+1}$ from the standard exponential distribution so that $U_{1,n}=_d E_1/\bar{E}_{n+1}$.
Furthermore
\begin{eqnarray*}
\log U_T(1/p)  &=&
\log T - {1 \over \alpha}log \left(1+ {p(n+1) \over k+1}\kappa (1+o(1))\right)\\
&&+
\log {\ell_U \left( (1/\bar{F}_W (T))[1+{p(n+1) \over k+1}\kappa (1+o(1))]^{-1}  \right) \over \ell_U \left( 1/\bar{F}_W (T)  \right)}.
\end{eqnarray*}
Now it follows from Propositions 1(a) and  2(a) that
$$
R_{k,n}^{\hat{\alpha}^T_{k,n}}  \to {1+\kappa}^{-1} \mbox{ and }
{(n+1)\hat D _T \over k+1} \to \kappa^{-1},
$$
so that, with
\begin{eqnarray*}
h_{\rho}\left([1+{\kappa \over k+1}E_1 (1+o_p(1))]^{-1}\right)
&=& -{\kappa \over k+1}E_1 (1+o_p(1)), \\
h_{\rho}\left([1+{p(n+1) \over k+1}\kappa (1+o_p(1))]^{-1}\right)
&=& -\kappa {p(n+1) \over k+1} (1+o_p(1)),
\end{eqnarray*}
\begin{eqnarray*}
\log \hat{q}^T_{p,k,n}-\log Q(1-p) &=&
\left(
-{1 \over \alpha}{\kappa \over k+1}E_1
+ {\kappa \over \hat{\alpha}^T_{k,n}}\left({1 \over k+1}-{p(n+1) \over k+1} \right)
\right. \\
&& \left.
+{\kappa  \over \alpha}{p(n+1) \over k+1}
- \kappa b(1/\bar{F}_W (T))\left({E_1 \over k+1}- {p(n+1) \over k+1}\right)
\right)(1+o_p(1)),
\end{eqnarray*}
from which (a) follows.

\vspace{1cm} \noindent
Next, in the case $k/(nD_T) \to \infty$, starting from expression (23) and Proposition 1(b), we obtain
$$
\log q^T_{p,k,n}=
\log X_{n-k,n} + {1 \over \alpha^T_{k,n}}\log \left(1+ {(n+1)\hat{D}_{T}\over k+1} \right) -
{1 \over \alpha^T_{k,n}}\log \left( {(n+1)p \over k+1}+{(n+1)\hat{D}_T \over k+1}\right),
$$
where under the given conditions
\begin{eqnarray*}
{(n+1)\hat{D}^T_{k,n}\over k+1} &=& {R_{k,n}^{\hat{\alpha}^T_{k,n}}-{1 \over k+1} \over 1- R^{\hat{\alpha}^T_{k,n}}_{k,n}} \\
&=&
{1 \over k+1}\left( e^{-[E_{k,k} -\log (k+1)]}(1+o_p(1)) -1 \right)\\
&=& {1 \over k+1}\left( E\,(1+o_p(1)) -1 \right),
\end{eqnarray*}
since, using Theorem 1(b) and Proposition 2(b), and the fact that $\exp (-[E_{k,k}- \log (k+1)])$ asymptotically is standard exponentially distributed,
\begin{eqnarray*}
(k+1)R_{k,n}^{\hat{\alpha}^T_{k,n}} &=&
(k+1)^{-({\hat{\alpha}^T_{k,n} \over \alpha}-1)}
\exp \left( -{\hat{\alpha}^T_{k,n} \over \alpha}(E_{k,k}-\log (k+1))\right)\\
&& \hspace{0.2cm} \times \left({1+(n+1)D_T \over 1+{(n+1)D_T \over k+1}} \right)^{\hat{\alpha}^T_{k,n}/\alpha} \\
&& \hspace{0.2cm}
\times \exp \left( - \hat{\alpha}^T_{k,n}b(n/k)\left[{(n+1)D_T \over k+1}+
h_{\rho}([k^{-1}+ {(n+1)D_T \over (k+1)}]^{-1}) \right]  \right) \\
&=& E\,(1+o_p(1)) -1.
\end{eqnarray*}
Furthermore,
\begin{eqnarray*}
\log X_{n-k,n}-\log Q(1-p) &=&
-{1 \over \alpha}\log {(n+1)U_{k+1,n} \over k+1} +
{1 \over \alpha}\log {(n+1)p \over k+1}\\
&& -\log (1+{(n+1)D_T \over k+1 }{k+1 \over (n+1)U_{k+1,n}} )
+  {1 \over \alpha}\log \left( 1+ {D_T \over p}\right)\\
&& + \log \ell_U \left(C_T {n+1 \over k+1} {k+1 \over (n+1)U_{k+1,n}}[1+{(n+1)D_T \over k+1} ]^{-1} \right)\\
&& - \log \ell_U \left(C_T \, p^{-1}[1+{D_T \over p} ]^{-1} \right).
\end{eqnarray*}
Finally since under the given assumptions $(n+1)U_{k+1,n}/(k+1) \to_p 1$, $D_T/p \to 0$, and
$$
\log {\ell_U \left(C_T {n+1 \over k+1} {k+1 \over (n+1)U_{k+1,n}}[1+{(n+1)D_T \over k+1} ]^{-1} \right)\over \ell_U \left(C_T \, p^{-1}[1+{D_T \over p} ]^{-1} \right)} = O(b(n/k)),
$$
we find that
\begin{eqnarray*}
\log \hat{q}^T_{p,k,n}-\log Q(1-p) &=&
\left\{ - \left({1 \over \hat{\alpha}^T_{k,n}}- {1 \over \alpha} \right)
\log ({(n+1)p \over k+1}) \right.\\&&\left.
 \hspace{0.7cm}-
{1 \over \hat{\alpha}^T_{k,n}} \log \left(1+ {E-1 \over (n+1)p} \right)
\right\}
(1+o_p(1)),
\end{eqnarray*}
from which the result follows with Theorem 1(b).
\end{proof}

\vspace{0.2cm} \noindent
\textbf{Proof of Theorem 3.} The first statement in (a) and (b) follow readily from Propositions 1(b) and 2(b).
\\\\
In order to prove (c) we first derive the asymptotic distribution of $\sqrt{k}\left( E_{k,n} (1/H_{k,n})- 1/2 \right)$.
Note that
\begin{eqnarray*}
\sqrt{k}\left( E_{k,n} (1/H_{k,n})- 1/2 \right) &=& \sqrt{k}\left( E_{k,n} (1/H_{k,n})- {1 \over 1+ {1 \over \alpha H_{k,n}} }\right)\\
 && \hspace{1cm} + \sqrt{k} \left({1 \over 1+ {1 \over \alpha H_{k,n} }} -1/2 \right)  \\
& =& \mathbb{E}_{k,n}(1/H_{k,n}) + {\Gamma _{k,n} \over 2(H_{k,n}+ \alpha^{-1})}
\end{eqnarray*}
with
\begin{eqnarray*}
\Gamma _{k,n} &=& \sqrt{k} (H_{k,n}-\alpha^{-1} ), \\
\mathbb{E}_{k,n}(s) &=& \sqrt{k}\left( E_{k,n} (s)- {1 \over 1+ s/\alpha }\right), \;\; s >0.
\end{eqnarray*}
In Theorem A.1 in Beirlant {\it et al.} (2009) \cite{Beirlant09} it is stated that $(\Gamma_{k,n}, \mathbb{E}_{k,n})$ converges weakly in the space $\mathbb{R} \times  \mathcal{C} [0, s_0] $ with $s_0 > 0$ and
$ \mathcal{C} [a,b]$ the Banach space of continuous functions $f: [a,b] \to \mathbb{R}$ equipped with the topolgy of uniform convergence. The limit process $(\Gamma, \mathbb{E})$ is a Gaussian process with
\begin{eqnarray*}
Var (\Gamma) &=& \alpha^{-2},\;\; Cov (\mathbb{E}(s_1),\mathbb{E}(s_2)) = {s_1s_2/\alpha^{2} \over (1+s_1/\alpha +s_2/\alpha)(1+s_1/\alpha)(1+s_2/\alpha)} \\
Cov (\Gamma, \mathbb{E}(s)) &=& {-s/\alpha^{2} \over (1+s/\alpha)^2} .
\end{eqnarray*}
From $H_{k,n} = \alpha^{-1} +o_p(1)$ we have for $s_0 > \alpha$ that $\mathbb{P} (s_0 \geq 1/H_{k,n}\geq 0) \to 1$ as $n \to \infty$, and thus by the continuous mapping theorem
$(\Gamma_{k,n},\mathbb{E}_{k,n}(1/H_{k,n}))$ converges weakly to $(\Gamma,\mathbb{E}(\alpha))$.
From this it then follows that
$$
 \mathbb{E}_{k,n}(1/H_{k,n}) + {\Gamma _{k,n} \over 2(H_{k,n}+ \alpha^{-1})} \to \mathcal{N}\left(0, 1/48 \right).
$$
Then also
$$
1- E_{k,n}(1/H_{k,n}) \to_p 1/2,
$$
from which the result (c) follows.
\\\\
To prove the second statement in (a), we obtain from (15) that
\begin{eqnarray*}
E_{k,n}(1/H_{k,n}) &=& {1 \over k}\sum_{j=1}^k \left( {Q_T(1-U_{k+1,n}) \over Q_T(1-U_{j,n}) }\right)^{1/H_{k,n}} \\
&=& {1 \over k}\sum_{j=1}^k \left( {1+  {U_{j,n}\over U_{k+1,n} }{U_{k+1,n}\over D_T }  \over 1+ {U_{k+1,n}\over D_T }}\right)^{1/(\alpha H_{k,n})} \\
& & \hspace{2cm} \times \left({ \ell_U({1 \over \bar{F}_W(T)}[1+ {U_{k+1,n}\over D_T } ]^{-1} )  \over  \ell_U({1 \over \bar{F}_W(T)}[1+ {U_{j,n}\over U_{k+1,n} }{U_{k+1,n}\over D_T } ]^{-1} )  } \right)^{1/H_{k,n}}.
\end{eqnarray*}
Setting $V_{j,k} := U_{j,n}/U_{k+1,n}$ ($1 \leq j \leq k$) we have that $V_{1,k}, \ldots, V_{k,k}$ are distributed as the order statistics of a uniform ($0,1$) sample of size $k$. Moreover for $U_{k+1,n}/D_T =: \kappa_{k,n}$ we have
$\kappa_{k,n} \to_p \kappa$.   Following Proposition 1(a) it then follows that
$$
H_{k,n} \to_p {1 \over \alpha}\left( 1- {1 \over \kappa}\log (1+\kappa) \right),
$$
\begin{eqnarray*}
&& \hspace{-1.5cm}
 \left( {1+  {U_{j,n}\over U_{k+1,n} }{U_{k+1,n}\over D_T }  \over 1+ {U_{k+1,n}\over D_T }}\right)^{1/(\alpha H_{k,n})} =
\left( {1+ \kappa_{k,n} V_{j,k} \over 1+ \kappa_{k,n} }\right)^{\kappa/(\kappa - \log (1+\kappa))(1+o(1))}, \\
&& \hspace{-1.5cm}
\left({ \ell_U({1 \over \bar{F}_W(T)}[1+ {U_{j,n}\over U_{k+1,n} }{U_{k+1,n}\over D_T } ]^{-1} )  \over  \ell_U({1 \over \bar{F}_W(T)}[1+ {U_{k+1,n}\over D_T } ]^{-1} )  } \right)^{-1/H_{k,n}} \\
&\sim &
1- {1 \over H_{k,n}} b(1/\bar{F}_W(T)) \left(h_{\rho}([1+\kappa V_{j,k}]^{-1}) -  h_{\rho}([1+\kappa]^{-1})\right).
\end{eqnarray*}
Since
\begin{eqnarray*}
 && \hspace{-2cm} {1 \over k}\sum_{j=1}^k \left( {1+ \kappa V_{j,k} \over 1+ \kappa }\right)^{\kappa /(\kappa -\log (1+\kappa))}\\  &\to_p &  \int_0^1 \left({1+\kappa u \over 1+\kappa}\right)^{\kappa /(\kappa -\log (1+\kappa))}du \\
 &  & \hspace{1cm}
= {(1+\kappa)(\kappa -\log (1+\kappa)) \over \kappa( 2\kappa -\log (1+\kappa))}
\left(1-(1+\kappa)^{-{2\kappa-\log (1+\kappa) \over \kappa -\log (1+\kappa)}} \right) \in (0, 1/2),
\end{eqnarray*}
and
$$
{1 \over k}\sum_{j=1}^k  (1+\kappa_{k,n}V_{j,k}) ^{\kappa /(\kappa -\log (1+\kappa))} \left(h_{\rho}([1+\kappa_{k,n}V_{j,k}]^{-1}) -  h_{\rho}([1+\kappa_{k,n}]^{-1})\right) =O_p(1),
$$
the result follows.

\end{document}